\documentclass[10pt]{article}

\usepackage{graphicx, amssymb, latexsym, amsfonts, amsmath, lscape, amscd, multirow, multicol,
amsthm, color, epsfig, mathrsfs, tikz, enumerate, soul, subcaption, scrextend, todonotes, xstring, hyperref, comment}

\newtheorem{theorem}{Theorem}
\newtheorem{conjecture}[theorem]{Conjecture}

\newtheorem{corollary}[theorem]{Corollary}

\newtheorem{lemma}[theorem]{Lemma}
\newtheorem{proposition}[theorem]{Proposition}
\newtheorem{question}{Question}
\newtheorem{problem}{Problem}

\newtheorem{observation}[theorem]{Observation}
\newtheorem{remark}[theorem]{Remark}

\theoremstyle{definition}
\newtheorem{claim}{Claim}

\let\oldenumerate\enumerate
\renewcommand{\enumerate}{
  \oldenumerate
  \setlength{\itemsep}{0pt}
  \setlength{\parskip}{0pt}
  \setlength{\parsep}{0pt}
}

\newenvironment{unnumbered}[1]{\trivlist \item [\hskip \labelsep {\bf
#1}]\ignorespaces\it}{\endtrivlist}

\newcommand\DELETE[1]{}

\newcommand{\1}{\vspace{0.1cm}}
\newcommand{\2}{\vspace{0.2cm}}

\newenvironment{proofofclaim}{\vspace{-6mm} \paragraph{\normalfont \textit{Proof of Claim.}} \hspace{-2.5mm}}{\hfill$\blacklozenge$ \medskip}

\newcommand{\tl}[1]{\textcolor{red}{#1}}
\newcommand{\ID}{\gamma^{{\rm ID}}}

\newcommand{\diam}{{\rm diam}}
\newcommand{\smallo}{{\rm o}}

\title{Identifying codes in graphs of given maximum degree: Characterizing trees\footnote{The first two authors were supported by the French government IDEX-ISITE initiative CAP 20-25 (ANR-16-IDEX-0001), the International Research Center "Innovation Transportation and Production Systems" of the I-SITE CAP 20-25, and the ANR project GRALMECO (ANR-21-CE48-0004). The research of Michael Henning was supported in part by the South African National Research Foundation (grant number 129265) and the University of Johannesburg. Tuomo Lehtil\"a's research was partially supported by Business Finland project 10769/31/2022 6GTNF and Academy of Finland grants 338797 and 358718.}~\footnote{Some results on the topic of this research and on a subclass of bipartite graphs have been announced in an extended abstract in the proceedings of the LAGOS 2023 conference~\cite{LAGOS}. We have since then been able to extend the research to all triangle-free graphs (in Part II), leading to the current presentation for Part I (on trees) in this paper.}}

\author{Dipayan Chakraborty\footnote{\noindent Université Clermont Auvergne, CNRS, Mines Saint-Étienne, Clermont Auvergne INP, LIMOS, 63000 Clermont-Ferrand, France.}~\footnote{\noindent Department of Mathematics and Applied Mathematics, University of Johannesburg, South Africa.}
\and Florent Foucaud\footnotemark[3]
\and Michael A. Henning\footnotemark[4]
\and Tuomo Lehtil\"a\footnote{University of Helsinki, Department of Computer Science, Helsinki, Finland.}~\footnote{Helsinki Institute for Information Technology (HIIT), Espoo, Finland.}~\footnote{University of Turku, Department of Mathematics and Statistics, Turku, Finland.}~~\footnotemark[3]
}

\date{\today}

\begin{document}

\maketitle

\begin{abstract}
An \emph{identifying code} of a closed-twin-free graph $G$ is a dominating set $S$ of vertices of $G$ such that any two vertices in $G$ have a distinct intersection between their closed neighborhoods and $S$. It was conjectured that there exists an absolute constant $c$ such that for every connected graph $G$ of order $n$ and maximum degree $\Delta$, the graph $G$ admits an identifying code of size at most $( \frac{\Delta-1}{\Delta} )n +c$. We provide significant support for this conjecture by exactly characterizing every tree requiring a positive constant $c$ together with the exact value of the constant. Hence, proving the conjecture for trees. For $\Delta=2$ (the graph is a path or a cycle), it is long known that $c=3/2$ suffices. For trees, for each $\Delta\ge 3$, we show that $c=1/\Delta\le 1/3$ suffices and that $c$ is required to have a positive value only for a finite number of trees. In particular, for $\Delta = 3$, there are 12 trees with a positive constant $c$ and, for each $\Delta \ge 4$, the only tree with positive constant $c$ is the $\Delta$-star. Our proof is based on induction and utilizes recent results from [F. Foucaud, T. Lehtil\"a. Revisiting and improving upper bounds for identifying codes. SIAM Journal on Discrete Mathematics, 2022]. We remark that there are infinitely many trees for which the bound is tight when $\Delta=3$; for every $\Delta\ge 4$, we construct an infinite family of trees of order $n$ with identification number very close to the bound, namely $\left( \frac{\Delta-1+\frac{1}{\Delta-2}}{\Delta+\frac{2}{\Delta-2}} \right) n > (\frac{\Delta-1}{\Delta} ) n -\frac{n}{\Delta^2}$. Furthermore, we also give a new tight upper bound for identification number on trees by showing that the sum of the domination and identification numbers of any tree $T$ is at most its number of vertices.
\end{abstract}


\section{Introduction}

A \emph{dominating set} of a graph $G$ is a set $S$ of vertices of $G$ such that every vertex in $G$ is dominated by a vertex in $S$, where a vertex \emph{dominates} itself and its neighbors. An \emph{identifying code} is a dominating set of a graph such that any two vertices are dominated by distinct subsets of vertices from the identifying code. Identifying codes were introduced in 1998 by Karpovsky, Chakrabarty and Levitin~\cite{karpovsky1998new}, motivated by fault-detection in multiprocessor networks. Numerous other applications of identifying codes have been discovered such as threat location in facilities using sensor networks~\cite{UTS04}, logical definability of graphs~\cite{PVV06} and canonical labeling of graphs for the graph isomorphism problem~\cite{B80}. Besides, since the 1960s and long before the introduction of identifying codes of graphs, many related concepts such as \emph{separating systems} or \emph{test covers} have been independently discovered and studied on hypergraphs: see  R\'enyi~\cite{renyi1961}, Bondy~\cite{bondy1972induced}, or Moret and Shapiro~\cite{MS85} for some early references. All of them put together form the general area of identification problems in graphs and other discrete structures, see the online bibliography~\cite{LocationBib} for over 500 papers on the topic.

Our goal is to upper-bound the smallest size of an identifying code in a graph of given maximum degree, motivated by a conjecture on this topic (see Conjecture~\ref{conj_G Delta_UB} below). In this paper, we prove the conjecture for trees, and characterize those trees that are extremal for the conjectured upper bound.

Let $G=(V(G),E(G))$ be a graph with vertex-set $V(G)$ and edge-set $E(G)$. Any subset $S \subset V(G)$ is called a \emph{vertex subset} of $G$. The \emph{(open) neighborhood} of a vertex $v$ of $G$ is the set $N_G(v)$ of all vertices of $G$ adjacent to $v$. The vertices of $G$ in $N_G(v)$ are also called the 
\emph{neighbors} of $v$ in $G$. Moreover, the set $N_G[v] = \{v\} \cup N_G(v)$ is called the \emph{closed neighborhood} of $v$. Vertices $u,v\in V(G)$ are called open (closed) twins in $G$ if and only if they have the same open (closed) neighborhood. Graphs with no open or closed twins are called \textit{twin-free}.

Formally speaking, an \emph{identifying code} $C$ of a graph $G$ is a vertex subset of $G$ that (i) \emph{dominates} each vertex $v$ of $G$ (that is, either $v \in C$ or  $v$ has a neighbor in $C$) and (ii) \emph{separates} each pair $u,v$ of distinct vertices of $G$ (that is, there is a vertex of $G$ in $C$ that belongs to exactly one of the two closed neighborhoods $N[u]$, $N[v]$). Note that graph $G$ admits an identifying code only when no two vertices of $G$ are closed twins. Hence, we say that graphs with no closed twins are \emph{identifiable}, that is, they admit an identifying code (for example, the whole vertex set). A vertex subset of a graph $G$ satisfying the property (i) is called a \emph{dominating set} of $G$ and a subset satisfying property (ii) is called a \textit{separating set} of $G$. It is natural to ask for a minimum-size identifying code of an identifiable graph $G$. The size, denoted by $\ID(G)$, of such a minimum-size identifying code of $G$ is called the \emph{identification number} of $G$. The smallest size of a dominating set of $G$ is denoted by $\gamma(G)$.

A natural question that arises in the study of identifying codes (like for the usual dominating sets) is the one of extremal values for the identification number: how large can it be, with respect to some relevant graph parameters? When only the order $n$ of the graph is considered, it is known that the identification number of an identifiable graph with at least one edge lies between $\log_2(n+1)$~\cite{karpovsky1998new} and $n-1$~\cite{GM07}; both values are tight and the extremal examples have been characterized in~\cite{M06} and~\cite{FGKNPV11}, respectively.

A thorough treatise on domination in graphs can be found in \cite{HaHeHe-20,HaHeHe-21,HaHeHe-23,HeYe-book}. Bounds on domination numbers for graphs with restrictions on their degree parameters are a natural and important line of research. In 1996, Reed~\cite{reed_1996} proved the influential result that if $G$ is a connected cubic graph of order~$n$, then $\gamma(G) \le \frac{3}{8}n$. In 2009, this bound was improved to $\gamma(G) \le \frac{5}{14}n$, if we exclude the two non-planar cubic graphs of order~$8$~\cite{KOSTOCHKA20091142}. A study of independent domination in graphs with bounded maximum degree has received considerable attention in the literature. We refer to the breakthrough paper~\cite{IDom-bound1}, as well as the papers in~\cite{IDom-bound2,IDom-bound3,IDom-bound4}. Another classical result is that the total domination number of a cubic graph is at most one-half its order~\cite{Tdom-bound1}, and the remarkable result that the total domination number of a $4$-regular graph is at most three-sevenths its order was proved by an interplay with the notion of transversals in hypergraphs~\cite{Tdom-bound2}. A detailed discussion on upper bounds on domination parameters in graphs in terms of their order and minimum and maximum degree, as well as bounds with specific structural restrictions imposed, can be found in~\cite[Chapters 6, 7, 10]{HaHeHe-23}.

With respect to the identification number, it was observed in~\cite{FGKNPV11} that when the maximum degree $\Delta$ of the graph $G$ is small enough with respect to the order $n$ of the graph, the ($n-1$)-upper bound can be significantly improved (for connected graphs) to $n-\frac{n}{\Theta(\Delta^5)}$. The latter was thereafter subsequently reduced to $n-\frac{n}{\Theta(\Delta^3)}$ in~\cite{foucaud2012degree}. This raises the question of what is the largest possible identification number of a connected identifiable graph of order $n$ and maximum degree $\Delta$. Towards this problem, the following conjecture was posed; the goal of this paper is to study this conjecture.

\begin{conjecture}[{\cite[Conjecture 1]{foucaud2012size}}] \label{conj_G Delta_UB}
There exists a constant $c$ such that for every connected identifiable graph of order $n\ge 2$ and of maximum degree $\Delta\ge 2$,
\[
\ID(G) \le \left( \frac{\Delta - 1}{\Delta} \right) n  + c.
\]
\end{conjecture}

Note that for $\Delta\le 1$, the only connected identifiable graph is the one-vertex graph. From the known results in the literature, the above conjecture holds for $\Delta=2$ (that is, for paths and cycles) with $c=3/2$ (see Theorem~\ref{lem_paths_cycles}). Hence, for the rest of this manuscript, we assume that $\Delta \ge 3$.

If true, Conjecture~\ref{conj_G Delta_UB} would be tight: for any value of $\Delta\ge 3$, there are arbitrarily large graphs of order $n$ and maximum degree $\Delta$ with identification number $( \frac{\Delta - 1}{\Delta} ) n$, see~\cite{foucaud2012degree,sierpinski}. A bound of the form $n-\frac{n}{103(\Delta+1)^3}$~\cite{foucaud2012degree} proved using probabilistic arguments is the best known general result towards Conjecture~\ref{conj_G Delta_UB} (for the sake of comparison, the conjectured bound can be rewritten as $n-\frac{n}{\Delta}+c$). It is reduced to $n-\frac{n}{103\Delta}$ for graphs with no \emph{forced} vertices (vertices that are the only difference between the closed neighborhoods of two vertices of the graph), and to $n-\frac{n}{f(k)\Delta}$ for graphs of clique number $k$~\cite{foucaud2012degree}. For triangle-free graphs, this was improved to $n-\frac{n}{\Delta+ \smallo(\Delta)}$ in~\cite{foucaud2012size}, and to smaller bounds for subclasses of triangle-free graphs, such as $n-\frac{n}{\Delta+9}$ for bipartite graphs and $n-\frac{n}{3\Delta/(\ln\Delta-1)}$ for triangle-free graphs without (open) twins. The latter result implies that Conjecture~\ref{conj_G Delta_UB} holds for triangle-free graphs without any open twins, whenever $\Delta\ge 55$ (because then, $3\Delta/(\ln\Delta-1)\le \Delta$).  Conjecture~\ref{conj_G Delta_UB} is also known to hold for line graphs of graphs of average degree at least~5~\cite[Corollary 21]{FGNPV13} as well as graphs which have girth at least $5$, minimum degree at least $2$ and maximum degree at least $4$ \cite{balbuena2015locating}. Moreover, it holds for bipartite graphs without (open) twins by~\cite{FL22}. Furthermore, the conjecture holds in many cases for some graph products such as Cartesian and direct products~\cite{goddard2013id, junnila2019conjecture, rall2014identifying}. See also the book chapter~\cite{lobstein2020locating}, where Conjecture~\ref{conj_G Delta_UB} is presented.

Until this work, the conjecture remained open even for trees, and one of the challenges of proving it on trees was to allow open-twins of degree~1, which are present in many trees (note that for any set of mutual open-twins, one needs all of them but one in any identifying code). We will also see that almost all extremal trees for Conjecture~\ref{conj_G Delta_UB} (those requiring $c>0$) have twins.
Moreover, it is known that when a tree $T$ of order at least~3 (except the path $P_4$) has no open twins, it satisfies $\ID(T)\le \frac{2}{3}n$~\cite[Corollary 8]{FL22} (this even holds for bipartite graphs). This implies the conjecture for $\Delta\ge 3$ for twin-free bipartite graphs, and clearly shows that the presence of open twins is the main difficulty in proving Conjecture~\ref{conj_G Delta_UB} for trees.

Different aspects of identifying codes in trees have previously been studied in the literature, see for example~\cite{AUGER20101372,BCHL2004,bertrand20051,BCMMS07,FL22,gimbel2001location,RRM19} for some examples.

\medskip

\noindent \textbf{Our work.} For graphs with maximum degree $\Delta = 2$ (that is, for paths and cycles), Conjecture~\ref{conj_G Delta_UB} is settled in~\cite{BCHL2004} and~\cite{GMS2006} with the constant $c = \frac{3}{2}$. We therefore investigate Conjecture~\ref{conj_G Delta_UB} when $\Delta \ge 3$.

Our aim is to prove Conjecture~\ref{conj_G Delta_UB} for all triangle-free graphs. To ease and shorten the presentation of our proof, we split it into two independent papers. In the current paper, we will focus on trees, and prove Conjecture~\ref{conj_G Delta_UB} (with $c= \frac{1}{\Delta}\le \frac{1}{3}$) for all trees of maximum degree at least~3 (the conjecture already holds with $c=1$ for all identifiable trees of maximum degree~2, that is, for paths of order at least~3, due to Bertrand et al\tl{.}~\cite{BCHL2004}).

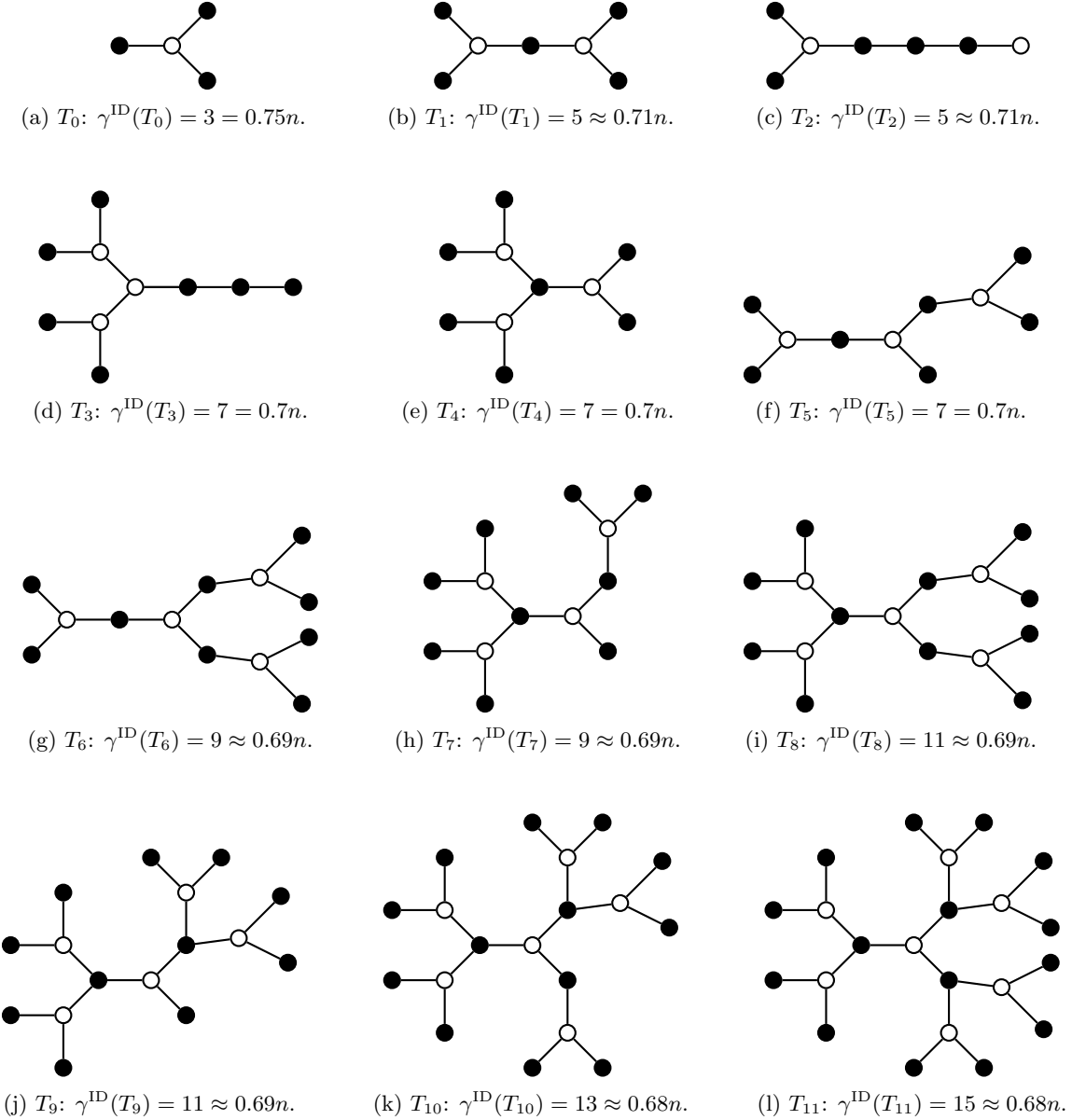
\begin{figure}[!htb]
\centering
\begin{subfigure}[t]{0.3\textwidth}
\centering
\begin{tikzpicture}[
blacknode/.style={circle, draw=black!, fill=black!, thick},
whitenode/.style={circle, draw=black!, fill=white!, thick},
scale=0.5]
\tiny
\node[blacknode] (0) at (0,0) {};
\node[whitenode] (1) at (1.5,0) {};
\node[blacknode] (2) at (2.5,1) {};
\node[blacknode] (3) at (2.5,-1) {};
\draw[-, thick] (0) -- (1);
\draw[-, thick] (1) -- (2);
\draw[-, thick] (1) -- (3);
\end{tikzpicture}
\caption{$T_0$: $\ID(T_0) = 3 = 0.75n$.}
\end{subfigure}
\hspace{2mm}
\begin{subfigure}[t]{0.3\textwidth}
\centering
\begin{tikzpicture}[
blacknode/.style={circle, draw=black!, fill=black!, thick},
whitenode/.style={circle, draw=black!, fill=white!, thick},
scale=0.5]
\tiny
\node[blacknode] (0) at (0,0) {};
\node[whitenode] (1) at (1.5,0) {};
\node[blacknode] (2) at (2.5,1) {};
\node[blacknode] (3) at (2.5,-1) {};
\node[whitenode] (1') at (-1.5,0) {};
\node[blacknode] (2') at (-2.5,1) {};
\node[blacknode] (3') at (-2.5,-1) {};
\draw[-, thick] (0) -- (1);
\draw[-, thick] (1) -- (2);
\draw[-, thick] (1) -- (3);
\draw[-, thick] (0) -- (1');
\draw[-, thick] (1') -- (2');
\draw[-, thick] (1') -- (3');
\end{tikzpicture}
\caption{$T_1$: $\ID(T_1) = 5 \approx 0.71n$.}
\end{subfigure}
\hspace{2mm}
\begin{subfigure}[t]{0.3\textwidth}
\centering
\begin{tikzpicture}[
blacknode/.style={circle, draw=black!, fill=black!, thick},
whitenode/.style={circle, draw=black!, fill=white!, thick},
scale=0.5]
\tiny
\node[whitenode] (1) at (0,0) {};
\node[blacknode] (2) at (-1,1) {};
\node[blacknode] (3) at (-1,-1) {};
\node[blacknode] (4) at (1.5,0) {};
\node[blacknode] (5) at (3,0) {};
\node[blacknode] (6) at (4.5,0) {};
\node[whitenode] (7) at (6,0) {};
\draw[-, thick, black!] (1) -- (2);
\draw[-, thick, black!] (1) -- (3);
\draw[-, thick, black!] (1) -- (4);
\draw[-, thick] (5) -- (4);
\draw[-, thick] (4) -- (6);
\draw[-, thick] (6) -- (7);
\end{tikzpicture}
\caption{$T_2$: $\ID(T_2) = 5 \approx 0.71n$.}
\end{subfigure}\vspace{8mm}
\hspace{2mm}
\begin{subfigure}[t]{0.3\textwidth}
\centering
\begin{tikzpicture}[
blacknode/.style={circle, draw=black!, fill=black!, thick},
whitenode/.style={circle, draw=black!, fill=white!, thick},
scale=0.5]
\tiny
\node[whitenode] (1) at (-1,1) {};
\node[blacknode] (2) at (-1,2.5) {};
\node[blacknode] (3) at (-2.5,1) {};
\node[whitenode] (4) at (0,0) {};
\node[whitenode] (8) at (-1,-1) {};
\node[blacknode] (9) at (-1,-2.5) {};
\node[blacknode] (10) at (-2.5,-1) {};
\node[blacknode] (5) at (1.5,0) {};
\node[blacknode] (6) at (3,0) {};
\node[blacknode] (7) at (4.5,0) {};
\draw[-, thick, black!] (1) -- (2);
\draw[-, thick, black!] (1) -- (3);
\draw[-, thick, black!] (1) -- (4);
\draw[-, thick, black!] (4) -- (8);
\draw[-, thick, black!] (8) -- (9);
\draw[-, thick, black!] (8) -- (10);
\draw[-, thick] (5) -- (4);
\draw[-, thick] (4) -- (6);
\draw[-, thick] (6) -- (7);
\end{tikzpicture}
\caption{$T_3$: $\ID(T_3) = 7 = 0.7n$.}
\end{subfigure}
\hspace{2mm}
\begin{subfigure}[t]{0.3\textwidth}
\centering
\begin{tikzpicture}[
blacknode/.style={circle, draw=black!, fill=black!, thick},
whitenode/.style={circle, draw=black!, fill=white!, thick},
scale=0.5]
\tiny
\node[blacknode] (0) at (0,0) {};
\node[whitenode] (1) at (1.5,0) {};
\node[blacknode] (2) at (2.5,1) {};
\node[blacknode] (3) at (2.5,-1) {};
\node[whitenode] (4) at (-1,-1) {};
\node[blacknode] (5) at (-1,-2.5) {};
\node[blacknode] (6) at (-2.6,-1) {};
\node[whitenode] (7) at (-1,1) {};
\node[blacknode] (8) at (-1,2.5) {};
\node[blacknode] (9) at (-2.6,1) {};
\draw[-, thick] (0) -- (1);
\draw[-, thick] (1) -- (2);
\draw[-, thick] (1) -- (3);
\draw[-, thick] (0) -- (4);
\draw[-, thick] (4) -- (5);
\draw[-, thick] (4) -- (6);
\draw[-, thick] (0) -- (7);
\draw[-, thick] (7) -- (8);
\draw[-, thick] (7) -- (9);
\end{tikzpicture}
\caption{$T_4$: $\ID(T_4) = 7 = 0.7n$.}
\end{subfigure}
\hspace{0mm}
\begin{subfigure}[t]{0.3\textwidth}
\centering
\begin{tikzpicture}[
blacknode/.style={circle, draw=black!, fill=black!, thick},
whitenode/.style={circle, draw=black!, fill=white!, thick},
scale=0.5]
\tiny
\node[blacknode] (0) at (0,0) {};
\node[whitenode] (1) at (1.5,0) {};
\node[blacknode] (2) at (2.5,1) {};
\node[blacknode] (3) at (2.5,-1) {};
\node[whitenode] (13) at (4,1.2) {};
\node[blacknode] (14) at (5.2,2.4) {};
\node[blacknode] (15) at (5.4,0.5) {};
\node[whitenode] (1') at (-1.5,0) {};
\node[blacknode] (2') at (-2.5,1) {};
\node[blacknode] (3') at (-2.5,-1) {};
\draw[-, thick] (0) -- (1);
\draw[-, thick] (1) -- (2);
\draw[-, thick] (1) -- (3);
\draw[-, thick] (2) -- (13);
\draw[-, thick] (13) -- (14);
\draw[-, thick] (13) -- (15);
\draw[-, thick] (0) -- (1');
\draw[-, thick] (1') -- (2');
\draw[-, thick] (1') -- (3');
\end{tikzpicture}
\caption{$T_5$: $\ID(T_5) = 7 = 0.7n$.}
\end{subfigure}\vspace{8mm}
\hspace{2mm}
\begin{subfigure}[t]{0.3\textwidth}
\centering
\begin{tikzpicture}[
blacknode/.style={circle, draw=black!, fill=black!, thick},
whitenode/.style={circle, draw=black!, fill=white!, thick},
scale=0.5]
\tiny
\node[blacknode] (0) at (0,0) {};
\node[whitenode] (1) at (1.5,0) {};
\node[blacknode] (2) at (2.5,1) {};
\node[blacknode] (3) at (2.5,-1) {};
\node[whitenode] (13) at (4,1.2) {};
\node[blacknode] (14) at (5.2,2.4) {};
\node[blacknode] (15) at (5.4,0.5) {};
\node[whitenode] (19) at (4,-1.2) {};
\node[blacknode] (20) at (5.2,-2.4) {};
\node[blacknode] (21) at (5.4,-0.5) {};
\node[whitenode] (1') at (-1.5,0) {};
\node[blacknode] (2') at (-2.5,1) {};
\node[blacknode] (3') at (-2.5,-1) {};
\draw[-, thick] (0) -- (1);
\draw[-, thick] (1) -- (2);
\draw[-, thick] (1) -- (3);
\draw[-, thick] (2) -- (13);
\draw[-, thick] (13) -- (14);
\draw[-, thick] (13) -- (15);
\draw[-, thick] (3) -- (19);
\draw[-, thick] (19) -- (20);
\draw[-, thick] (19) -- (21);
\draw[-, thick] (0) -- (1');
\draw[-, thick] (1') -- (2');
\draw[-, thick] (1') -- (3');
\end{tikzpicture}
\caption{$T_6$: $\ID(T_6) = 9 \approx 0.69n$.}
\end{subfigure}
\hspace{2mm}
\begin{subfigure}[t]{0.3\textwidth}
\centering
\begin{tikzpicture}[
blacknode/.style={circle, draw=black!, fill=black!, thick},
whitenode/.style={circle, draw=black!, fill=white!, thick},
scale=0.5]
\tiny
\node[blacknode] (0) at (0,0) {};
\node[whitenode] (1) at (1.5,0) {};
\node[blacknode] (2) at (2.5,1) {};
\node[blacknode] (3) at (2.5,-1) {};
\node[whitenode] (4) at (-1,-1) {};
\node[blacknode] (5) at (-1,-2.5) {};
\node[blacknode] (6) at (-2.5,-1) {};
\node[whitenode] (7) at (-1,1) {};
\node[blacknode] (8) at (-1,2.5) {};
\node[blacknode] (9) at (-2.5,1) {};
\node[whitenode] (10) at (2.5,2.5) {};
\node[blacknode] (11) at (3.5,3.5) {};
\node[blacknode] (12) at (1.5,3.5) {};
\draw[-, thick] (0) -- (1);
\draw[-, thick] (1) -- (2);
\draw[-, thick] (1) -- (3);
\draw[-, thick] (2) -- (10);
\draw[-, thick] (10) -- (11);
\draw[-, thick] (10) -- (12);
\draw[-, thick] (0) -- (4);
\draw[-, thick] (4) -- (5);
\draw[-, thick] (4) -- (6);
\draw[-, thick] (0) -- (7);
\draw[-, thick] (7) -- (8);
\draw[-, thick] (7) -- (9);
\end{tikzpicture}
\caption{$T_7$: $\ID(T_7) = 9 \approx 0.69n$.}
\end{subfigure}
\hspace{0mm}
\begin{subfigure}[t]{0.3\textwidth}
\centering
\begin{tikzpicture}[
blacknode/.style={circle, draw=black!, fill=black!, thick},
whitenode/.style={circle, draw=black!, fill=white!, thick},
scale=0.5]
\tiny
\node[blacknode] (0) at (0,0) {};
\node[whitenode] (1) at (1.5,0) {};
\node[blacknode] (2) at (2.5,1) {};
\node[blacknode] (3) at (2.5,-1) {};
\node[whitenode] (13) at (4,1.2) {};
\node[blacknode] (14) at (5.2,2.4) {};
\node[blacknode] (15) at (5.4,0.5) {};
\node[whitenode] (19) at (4,-1.2) {};
\node[blacknode] (20) at (5.2,-2.4) {};
\node[blacknode] (21) at (5.4,-0.5) {};
\node[whitenode] (4) at (-1,-1) {};
\node[blacknode] (5) at (-1,-2.5) {};
\node[blacknode] (6) at (-2.5,-1) {};
\node[whitenode] (7) at (-1,1) {};
\node[blacknode] (8) at (-1,2.5) {};
\node[blacknode] (9) at (-2.5,1) {};
\draw[-, thick] (0) -- (1);
\draw[-, thick] (1) -- (2);
\draw[-, thick] (1) -- (3);
\draw[-, thick] (2) -- (13);
\draw[-, thick] (13) -- (14);
\draw[-, thick] (13) -- (15);
\draw[-, thick] (3) -- (19);
\draw[-, thick] (19) -- (20);
\draw[-, thick] (19) -- (21);
\draw[-, thick] (0) -- (4);
\draw[-, thick] (4) -- (5);
\draw[-, thick] (4) -- (6);
\draw[-, thick] (0) -- (7);
\draw[-, thick] (7) -- (8);
\draw[-, thick] (7) -- (9);
\end{tikzpicture}
\caption{$T_8$: $\ID(T_8) = 11 \approx 0.69n$.}
\end{subfigure}\vspace{8mm}
\hspace{6mm}
\begin{subfigure}[t]{0.3\textwidth}
\centering
\begin{tikzpicture}[
blacknode/.style={circle, draw=black!, fill=black!, thick},
whitenode/.style={circle, draw=black!, fill=white!, thick},
scale=0.5]
\tiny
\node[blacknode] (0) at (0,0) {};
\node[whitenode] (1) at (1.5,0) {};
\node[blacknode] (2) at (2.5,1) {};
\node[blacknode] (3) at (2.5,-1) {};
\node[whitenode] (4) at (-1,-1) {};
\node[blacknode] (5) at (-1,-2.5) {};
\node[blacknode] (6) at (-2.5,-1) {};
\node[whitenode] (7) at (-1,1) {};
\node[blacknode] (8) at (-1,2.5) {};
\node[blacknode] (9) at (-2.5,1) {};
\node[whitenode] (10) at (2.5,2.5) {};
\node[blacknode] (11) at (3.5,3.5) {};
\node[blacknode] (12) at (1.5,3.5) {};
\node[whitenode] (13) at (4,1.2) {};
\node[blacknode] (14) at (5.2,2.4) {};
\node[blacknode] (15) at (5.4,0.5) {};
\draw[-, thick] (0) -- (1);
\draw[-, thick] (1) -- (2);
\draw[-, thick] (1) -- (3);
\draw[-, thick] (2) -- (10);
\draw[-, thick] (10) -- (11);
\draw[-, thick] (10) -- (12);
\draw[-, thick] (2) -- (13);
\draw[-, thick] (13) -- (14);
\draw[-, thick] (13) -- (15);
\draw[-, thick] (0) -- (4);
\draw[-, thick] (4) -- (5);
\draw[-, thick] (4) -- (6);
\draw[-, thick] (0) -- (7);
\draw[-, thick] (7) -- (8);
\draw[-, thick] (7) -- (9);

\end{tikzpicture}
\caption{$T_9$: $\ID(T_9) = 11 \approx 0.69n$.}
\end{subfigure}
\hspace{4mm}
\begin{subfigure}[t]{0.3\textwidth}
\centering
\begin{tikzpicture}[
blacknode/.style={circle, draw=black!, fill=black!, thick},
whitenode/.style={circle, draw=black!, fill=white!, thick},
scale=0.5]
\tiny
\node[blacknode] (0) at (0,0) {};
\node[whitenode] (1) at (1.5,0) {};
\node[blacknode] (2) at (2.5,1) {};
\node[blacknode] (3) at (2.5,-1) {};
\node[whitenode] (4) at (-1,-1) {};
\node[blacknode] (5) at (-1,-2.5) {};
\node[blacknode] (6) at (-2.5,-1) {};
\node[whitenode] (7) at (-1,1) {};
\node[blacknode] (8) at (-1,2.5) {};
\node[blacknode] (9) at (-2.5,1) {};
\node[whitenode] (10) at (2.5,2.5) {};
\node[blacknode] (11) at (3.5,3.5) {};
\node[blacknode] (12) at (1.5,3.5) {};
\node[whitenode] (13) at (4,1.2) {};
\node[blacknode] (14) at (5.2,2.4) {};
\node[blacknode] (15) at (5.4,0.5) {};
\node[whitenode] (16) at (2.5,-2.5) {};
\node[blacknode] (17) at (3.5,-3.5) {};
\node[blacknode] (18) at (1.5,-3.5) {};
\draw[-, thick] (0) -- (1);
\draw[-, thick] (1) -- (2);
\draw[-, thick] (1) -- (3);
\draw[-, thick] (2) -- (10);
\draw[-, thick] (10) -- (11);
\draw[-, thick] (10) -- (12);
\draw[-, thick] (2) -- (13);
\draw[-, thick] (13) -- (14);
\draw[-, thick] (13) -- (15);
\draw[-, thick] (3) -- (16);
\draw[-, thick] (16) -- (17);
\draw[-, thick] (16) -- (18);
\draw[-, thick] (0) -- (4);
\draw[-, thick] (4) -- (5);
\draw[-, thick] (4) -- (6);
\draw[-, thick] (0) -- (7);
\draw[-, thick] (7) -- (8);
\draw[-, thick] (7) -- (9);
\end{tikzpicture}
\caption{$T_{10}$: $\ID(T_{10}) = 13 \approx 0.68n$.}
\end{subfigure}
\hspace{4mm}
\begin{subfigure}[t]{0.3\textwidth}
\centering
\begin{tikzpicture}[
blacknode/.style={circle, draw=black!, fill=black!, thick},
whitenode/.style={circle, draw=black!, fill=white!, thick},
scale=0.5]
\tiny
\node[blacknode] (0) at (0,0) {};
\node[whitenode] (1) at (1.5,0) {};
\node[blacknode] (2) at (2.5,1) {};
\node[blacknode] (3) at (2.5,-1) {};
\node[whitenode] (4) at (-1,-1) {};
\node[blacknode] (5) at (-1,-2.5) {};
\node[blacknode] (6) at (-2.5,-1) {};
\node[whitenode] (7) at (-1,1) {};
\node[blacknode] (8) at (-1,2.5) {};
\node[blacknode] (9) at (-2.5,1) {};
\node[whitenode] (10) at (2.5,2.5) {};
\node[blacknode] (11) at (3.5,3.5) {};
\node[blacknode] (12) at (1.5,3.5) {};
\node[whitenode] (13) at (4,1.2) {};
\node[blacknode] (14) at (5.2,2.4) {};
\node[blacknode] (15) at (5.4,0.5) {};
\node[whitenode] (16) at (2.5,-2.5) {};
\node[blacknode] (17) at (3.5,-3.5) {};
\node[blacknode] (18) at (1.5,-3.5) {};
\node[whitenode] (19) at (4,-1.2) {};
\node[blacknode] (20) at (5.2,-2.4) {};
\node[blacknode] (21) at (5.4,-0.5) {};
\draw[-, thick] (0) -- (1);
\draw[-, thick] (1) -- (2);
\draw[-, thick] (1) -- (3);
\draw[-, thick] (2) -- (10);
\draw[-, thick] (10) -- (11);
\draw[-, thick] (10) -- (12);
\draw[-, thick] (2) -- (13);
\draw[-, thick] (13) -- (14);
\draw[-, thick] (13) -- (15);
\draw[-, thick] (3) -- (16);
\draw[-, thick] (16) -- (17);
\draw[-, thick] (16) -- (18);
\draw[-, thick] (3) -- (19);
\draw[-, thick] (19) -- (20);
\draw[-, thick] (19) -- (21);
\draw[-, thick] (0) -- (4);
\draw[-, thick] (4) -- (5);
\draw[-, thick] (4) -- (6);
\draw[-, thick] (0) -- (7);
\draw[-, thick] (7) -- (8);
\draw[-, thick] (7) -- (9);
\end{tikzpicture}
\caption{$T_{11}$: $\ID(T_{11}) = 15 \approx 0.68n$.}
\end{subfigure}

\caption{The family~$\mathcal{T}_{3}$ of trees of maximum degree~3 requiring $c>0$ in Conjecture~\ref{conj_G Delta_UB}. The set of black vertices in each figure constitutes an identifying code of the tree.}
\label{fig:trees}
\end{figure}

The main challenge for proving the conjecture, is the constant $c$ that could, in principle, be arbitrary. Thus, in order to prove it for trees, a large part of our proof is dedicated to analyzing those extremal trees that require $c>0$. In fact, for each $\Delta \ge 3$, we characterize the trees with maximum degree $\Delta$ for which $c>0$. The number of these trees is largest for $\Delta=3$, and in fact, this case is the hardest part of our proof. The characterization is given by the collection $\mathcal{T}_\Delta$ for $\Delta \ge 3$, where, for $\Delta=3$, $\mathcal{T}_\Delta$ is the set of 12 trees of maximum degree~$3$ and diameter at most~$6$ depicted in Figure~\ref{fig:trees}; and $\mathcal{T}_\Delta = \{K_{1,\Delta}\}$ for $\Delta \ge 4$.

In the companion paper~\cite{PaperPart2}, we will use the results from the current paper as the starting point of a proof of Conjecture~\ref{conj_G Delta_UB} for all triangle-free graphs (with the same list of trees as exceptional cases for which $c>0$ is necessary when $\Delta\ge 3$). Besides Theorem~\ref{the_main}, also our analysis about the exceptional trees from Figure~\ref{fig:trees} will be crucial for~\cite{PaperPart2}.

Note that, for maximum degree at least~$3$, all trees are identifiable. Hence, throughout the rest of the paper, we tacitly assume all our trees to be identifiable. Our main results are stated as follows.

\begin{theorem} \label{the_main}
Let $G$ be a tree of order $n$ and of maximum degree $\Delta \ge 3$. If $G$ is isomorphic to a tree in $\mathcal{T}_\Delta$, then, we have
\[
\ID(G) = \left( \frac{\Delta - 1}{\Delta} \right) n + \frac{1}{\Delta}.
\]
On the other hand, if $G$ is not isomorphic to any tree in the collection $\mathcal{T}_\Delta$, then we have
\[
\ID(G) \le \left( \frac{\Delta - 1}{\Delta} \right) n.
\]
\end{theorem}

We also determine the exact value of the identification number for the exceptional trees in $\mathcal{T}_\Delta$, as follows.

We have listed every tree requiring a positive constant $c$ for Conjecture \ref{conj_G Delta_UB} in Table \ref{TableConstVal}. In the companion paper \cite{PaperPart2}, we show that the only other triangle-free graphs which require a positive constant $c$ are odd cycles, and some small even-length cycles.

\begin{table}[h!]
\centering
\begin{tabular}{|c|c|c|c|}
\hline
Graph class & $\Delta$  & $c$  & Reference\\ \hline
\hline
$K_{1,\Delta}$ & $\Delta \ge 3$  & $1/\Delta$  &  Lemma \ref{lem_star}\\ \hline
$\mathcal{T}_{3}$ & 3  & $1/3$  & Theorem \ref{the_main} \\ \hline
Even paths & 2 & $1$  & \cite{BCHL2004} (Theorem \ref{lem_paths_cycles}) \\ \hline
Odd paths & 2  & $1/2$  & \cite{BCHL2004} (Theorem \ref{lem_paths_cycles})  \\ \hline
\end{tabular}
\caption{Trees requiring a positive constant $c$ for Conjecture \ref{conj_G Delta_UB}.}\label{TableConstVal}
\end{table}

We also show that Theorem~\ref{the_main} is tight for many trees with $c=0$, besides the list of exceptional trees mentioned above which require $c>0$. One tight example of trees with $c=0$ is the $2$-corona of a path~\cite{FL22}. We will see that other such tight examples exist. When $\Delta=3$, there are infinitely many examples for such trees. Furthermore, we give in Proposition~\ref{prop:BigConstruction} an infinite family of trees for any $\Delta\ge 4$ which have identification number quite close to the conjectured bound, namely $\frac{\Delta-1+\frac{1}{\Delta-2}}{\Delta+\frac{2}{\Delta-2}}n> (\frac{\Delta-1}{\Delta} ) n -\frac{n}{\Delta^2}$. In particular, as $\Delta$ increases, our construction gets closer and closer to the conjectured bound.

\medskip

\noindent \textbf{Structure of the paper.} Following the introduction in the current section, Section~\ref{sec:prelim} contains preliminary lemmas and results. In Section~\ref{sec:domBound}, we improve an upper bound from the literature, thereby relating the identifying code number of a tree with the usual domination number. In Section~\ref{sec:small_cases}, we study extremal trees leading to the proof of the first part of Theorem~\ref{the_main}. Section~\ref{sec:main} is then dedicated to the proof of Theorem~\ref{the_main}. In Section~\ref{sec:construct}, we propose some constructions to deal with the tightness of the conjectured bound for trees. We conclude in Section~\ref{sec:conclu}.

\section{Preliminaries}\label{sec:prelim}

In this section, we establish notations and mention some useful results from the literature.

\subsection{Definitions and notations}

For any vertex $v$ of $G$, the symbol $\deg_G(v)$ denoting the \emph{degree} of the vertex $v$ \emph{in} $G$ is the total number of neighbors of $v$ in $G$. A \emph{leaf} of a graph $G$ is a vertex of degree~$1$ in $G$. In the literature, a leaf is also known as a \emph{pendant vertex}. The (only) neighbor of a leaf $v$ in a graph $G$ is called the \emph{support vertex} of $v$ in $G$. The number of leaves and support vertices in $G$ are denoted by $\ell(G)$ and $s(G)$, respectively. Naturally, any vertex of a graph $G$ that is not a leaf of $G$ is usually referred to as a \emph{non-leaf} vertex of $G$. The length (or the number of edges) of a longest induced path in a graph $G$ is called the \emph{diameter} of $G$ and is denoted by $\diam(G)$.

On many occasions throughout this article, we shall have the need to look at a subgraph of a graph $G$ formed by deleting away some vertices or edges from $G$. To that end, given a graph $G$ and a set $S$ containing some vertices and edges of $G$, we define $G-S$ as the subgraph of $G$ obtained by deleting from $G$ all vertices (and edges incident with them) and edges of $G$ in $S$. For any positive integer $p$, let $[p]$ denote the set $\{1, 2, \ldots, p\}$ and for any two distinct integers $p$ and $q$ with $p < q$, let $[p,q]$ denote the set $\{p, p+1, \ldots , q\}$.

Given two sets $A$ and $B$, the set $A \ominus B = (A \setminus B) \cup (B \setminus A)$ is the \emph{symmetric difference} of $A$ and $B$. Now, let $C$ be a vertex subset of $G$. For a given vertex $v$ and a vertex subset $C$ of a graph $G$, the intersection $N_G[v] \cap C$, denoted by $C_G[v]$, is called the \emph{$C$-code of $v$ in $G$}; and each element of the $C$-code of $v$ in $G$ is called a \emph{$C$-codeword (of $v$ in $G$)}. Given a pair $u,v$ of distinct vertices of $G$ such that $C_G[u] \ominus C_G[v] \ne \emptyset$, a $C$-codeword of either $u$ or $v$ in $G$ in the non-empty symmetric difference $C_G[u] \ominus C_G[v]$ is called an \emph{separating $C$-codeword} for the pair $u,v$ in $G$.

One can check that the following remark holds.

\begin{remark}\label{rem_char_ID-code}
Let $G$ be an identifiable graph. A dominating set $C$ of $G$ is an identifying code of $G$ if and only if $C$ separates every pair $u,v$ of distinct vertices of $G$ such that $d_G(u,v) \le 2$.
\end{remark}

In light of the above Remark~\ref{rem_char_ID-code}, on all occasions throughout this paper where we need to prove a dominating set $C$ of an identifiable graph $G$ to be an identifying code of $G$, we simply show that $C$ separates all pairs of distinct vertices of $G$ with distance at most~2.

\subsection{Paths}

The domination number of a path $P_n$ on $n$ vertices is given by the following closed formula, noting that starting from the second vertex of the path, we add every third vertex on the path into our dominating set, together with the last vertex when $n$ is not divisible by~$3$.

\begin{observation}[{\cite[Observation~2.1]{HaHeHe-23}}]
\label{ob:path}
If $P_n$ is a path on $n$ vertices, then $\gamma(P_n) = \lceil \frac{n}{3} \rceil \le \frac{n+2}{3}$.
\end{observation}

It is known that Conjecture~\ref{conj_G Delta_UB} is settled with $c=1$ for all identifiable trees of maximum degree~$2$, that is, paths of order at least~$3$. This is due to the following result by Bertrand et al.~\cite{BCHL2004}. Observe that Conjecture \ref{conj_G Delta_UB} holds for even paths with $c=\frac{1}{2}$ and for odd paths with $c=1$.

\begin{theorem}[\cite{BCHL2004}] \label{lem_paths_cycles}
If $P_n$ is a path on $n$ vertices, then
\[
\ID(P_n) =
\begin{cases} \frac{n}{2}+\frac{1}{2}, &\text{if $n \ge 1$ is odd}, \1\\
								  \frac{n}{2}+1, &\text{if $n \ge 4$ is even}.
\end{cases}
\]
Therefore,
\[
\ID(P_n) = \Big \lfloor \frac{n}{2} \Big \rfloor+1 \le \left( \frac{\Delta - 1}{\Delta} \right) n+1.
\]
\end{theorem}

\subsection{Two useful lemmas from the literature}

We cite here the following two lemmas from~\cite{FL22}, that will be essential for the proof of our main result, Theorem~\ref{the_main}. The first lemma was originally proved only for trees in~\cite{haynes2006locating}.

\begin{lemma}[{\cite[Lemma 4]{FL22}}] \label{lem_FT_1}
If $G$ is a connected bipartite graph on $n \ge 4$ vertices with $s$ support vertices and not isomorphic to a path $P_4$, then
$\ID(G) \le n-s$.
\end{lemma}

\begin{lemma}[{\cite[Theorem 6]{FL22}}] \label{lem_FT_2}
If $G$ is a connected bipartite graph on $n \ge 3$ vertices and $\ell$ leaves, and with no twins of degree~2 or greater, then
$\ID(G) \le \frac{1}{2}(n+\ell)$.
\end{lemma}

Note that both Lemma~\ref{lem_FT_1} and Lemma~\ref{lem_FT_2} apply to trees, since trees are bipartite graphs in which the only possible twins are leaves sharing the same support vertex, and thus, they have degree~1.

\section{An upper bound relating the identification number of a tree with its domination number}\label{sec:domBound}

In this section, we introduce a simple new upper bound for identifying codes in trees, which improves the upper bound $\ID(T)\le n-s$ from Lemma~\ref{lem_FT_1} (for trees), since in any graph of order at least~3, the domination number is at least the number of support vertices. Note that the bound of Lemma~\ref{lem_FT_1} was initially proved in~\cite{haynes2006locating} for trees only.

\begin{theorem}\label{the_domBound}
If $T$ is a tree other than $P_4$ of order $n \ge 3$, then $\ID(T) \le n-\gamma(T)$.
\end{theorem}
\begin{proof}
We prove the claim by induction on the order of the tree. But first, we show that the claim holds for paths. By Observation~\ref{ob:path}, $\gamma(P_n) = \lceil \frac{n}{3} \rceil\le \frac{n+2}{3}$. Moreover, by Theorem \ref{lem_paths_cycles}, we have $\ID(P_n)= \frac{1}{2}(n+2)$ for even $n$ and $\ID(P_n)= \frac{1}{2}(n+1)$ for odd $n$. Furthermore, we have $\ID(P_n)+\gamma(P_n)\le \frac{n+2}{2}+\frac{n+2}{3}=\frac{5n+10}{6}\le n$, when $n\ge 10$ is even and $\ID(P_n)+\gamma(P_n)\le \frac{n+1}{2}+\frac{n+2}{3}=\frac{5n+7}{6}\le n$, when $n\ge 7$ is odd. Hence, the claim follows for even-length paths on at least ten vertices and for odd-length paths on at least seven vertices. Furthermore, we have $\ID(P_3)+\gamma(P_3)=2+1=3$, $\ID(P_5)+\gamma(P_5)=3+2=5$, $\ID(P_6)+\gamma(P_6)=4+2=6$ and $\ID(P_8)+\gamma(P_8)=5+3=8$. Hence, the claim follows for all paths.

For $n=3$, the only tree is $P_3$ which satisfies $\ID(P_3)=2$ and $\gamma(P_3)=1$, and so the bound holds. For $n=4$, if $T$ is not $P_4$, then it is the star $K_{1,3}$, with $\ID(K_{1,3})=3$ and $\gamma(K_{1,3})=1$. Hence, the claim is true for $n\le 4$.

Let us next assume that the claimed bound holds for every tree $T'$ other than $P_1,P_2,P_4$ that has order at most $n'$.
Suppose, to the contrary, that there are trees of order $n'+1$ for which the bound does not hold. Let $T$ be a tree of order $n=n'+1\ge 5$ such that $\ID(T)>n-\gamma(T)$. By the first part of the proof, the tree $T$ is not a path.

If all vertices of $T$ are either leaves or support vertices, then $\gamma(T)\le s(T)$ and $n-\gamma(T)\ge n-s(T)$ and the bound of the statement holds by Lemma \ref{lem_FT_1}, a contradiction. Thus, there exists a non-support, non-leaf vertex in $T$. Let us root the tree at such a non-leaf, non-support vertex $x$, and let us consider the support vertex $s$ which has the greatest distance to $x$. Notice that vertex $s$ is adjacent to exactly one non-leaf vertex due to its maximal distance to $x$.

Let us first assume that $s$ is adjacent to $t\ge 2$ leaves $l_1,l_2,\dots, l_t$. Consider the tree $T_s=T-\{s,l_1,l_2,\dots,l_t\}$. Since $T_s$ has $n-t-1$ vertices, by induction we have $\ID(T_s)\le n-t-1-\gamma(T_s)$ unless $T_s=P_4$, or $T_s$ contains at most two vertices. However, notice that $|V(T_s)| \ge 3$ on account of $x$ being a non-leaf and non-support vertex of $T$. Moreover, if $T_s = P_4$, then $s$ cannot be adjacent to a support vertex of $T_s$, or else, the vertex $x$ would be either a leaf or a support vertex of $T$, a contradiction to our assumption. Therefore the vertex $s$ must be adjacent to a leaf of $T_s$. In this case, we have $\gamma(T)=2=s(T)$. Since $\ID(T)\le n-s(T)$ by Lemma \ref{lem_FT_1}, the claim holds in this case. 
Therefore, we may assume by induction that there exists an identifying code $C_s$ in $T_s$ which has cardinality at most $n-t-1-\gamma(T_s)$, with $D_s$ a minimum-size dominating set in $T_s$. Furthermore, set $C_s\cup\{l_1,\dots,l_t\}$ is an identifying code of $T$, and set $D_s\cup\{s\}$ is a dominating set of $T$. Thus, $\ID(T)\le n-t-1-\gamma(T_s)+t\le n-\gamma(T)$, as claimed.

Let us assume next that support vertex $s$ has degree~2, and denote by $u$ the non-leaf adjacent to $s$ (possibly, $u=x$) and by $l$ the leaf adjacent to $s$. If vertex $u$ also has degree~2, then we consider the tree $T_u=T-\{l,s,u\}$. Since $x$ is a non-leaf and non-support vertex of $T$, we must have $|V(T_u)| \ge 2$. However, if $|V(T_u)| = 2$, then $T$ is a path, a contradiction. 
Moreover, if $T_u\cong P_4$, then $T$ is path if $u$ is adjacent to a leaf of $T_u$, again a contradiction. Hence, $T$ has three support vertices and $\ID(T)\le n-3$ (consider the identifying code formed by all non-leaf vertices of $T$) while $\gamma(T)=3$. Hence, by induction, we may assume that $T_u$ has an optimal identifying code $C_u$ with cardinality $|C_u|=\ID(T_u)\le n-3-\gamma(T_u)$ and we consider a minimum-size dominating set $D_u$ of $T_u$. Notice that $D_u\cup\{s\}$ is a dominating set in $T$ and either $C_u\cup\{s,u\}$ or $C_u\cup\{l,u\}$ is an identifying code in $T$. Hence,  $\ID(T)\le n-3-\gamma(T_u)+2\le n-\gamma(T)$ as claimed.

Therefore, we may assume from now on that $\deg_T(u)\ge 3$.

Let us next consider the case where $u$ is a support vertex with adjacent leaf $l_u$. Furthermore, let us consider the tree $T_s=T-\{s,l\}$. Notice that $T_s$ contains at least two vertices, and if $T_s$ contains only two vertices, then $T$ is a path, a contradiction. Moreover, if $T_s\cong P_4$, then $u$ is a support vertex of $P_4$ and $\ID(T)=\gamma(T)=3=n-\gamma(T)$ (consider the set of non-leaf vertices as an identifying code of $T$). Hence, we may apply induction to $T_s$ and thus, there exists an optimal identifying code $C_s$ of $T_s$ of cardinality $|C_s|=\ID(T_s)\le n-2-\gamma(T_s)$. Also consider a minimum-size dominating set $D_s$ in $T_s$. Assume first that $u\in C_s$. Now set $C_s\cup\{l\}$ is an identifying code in $T$ and $D_s\cup\{s\}$ is a dominating set in $T$. Thus, $\ID(T)\le |C_s|+1\le n-2-\gamma(T_s)+1\le n-\gamma(T)$ and we are done. Moreover, if $u\not\in C_s$, then $l_u\in C_s$. Furthermore, in this case set $C=(C_s\cup\{s,u\})\setminus\{l_u\}$ is an identifying code in $T$. Indeed, since $C_s$ is an identifying code in $T_s$ and $N_{T_s}[l_u]\cap C_s=\{l_u\}$, we have $|N_{T_s}[u]\cap C_s|\ge 2$. Hence, we have $|N_T[u]\cap C|\ge 3$, $N_{T}[l_u]\cap C=\{u\}$, $N_{T}[s]\cap C=\{u,s\}$ and  $N_{T}[l]\cap C=\{s\}$, confirming that $C$ is an identifying code of $T$. Moreover, we have $\ID(T)\le n-\gamma(T)$ by the same arguments as in the case where $u\in C_s$.

Hence, we may assume that $u$ is not a support vertex.

Since we assumed that $s$ is a support vertex with the greatest distance to $x$, there is at most one non-support vertex adjacent to $u$ (on the path from $u$ to $x$). Moreover, if there exists a support vertex $s'\in N(u)$ with $\deg(s')\ge 3$, then we could have considered $s'$ instead of $s$ in an earlier argument considering the case $\deg(s)\ge 3$. Thus, we may assume that every support neighbor of $u$ has degree~$2$. Let us denote the support vertices adjacent to $u$ by $\{s_1,\dots, s_h\}$ (with $s=s_1$) and the leaf adjacent to $s_i$ by $l_i$ for each $1\le i \le h$. Observe that the tree $SS_h=T[\{s_1,\dots,s_h,l_1,\dots,l_h,u\tl{\}}]$ is a subdivided star. Moreover, it has domination number $h$ and identification number $h+1$. Indeed, we may take the center and every support vertex as our identifying code $C_{h}$. Hence, the subdivided star $SS_h$ satisfies the claimed upper bound: $\ID(SS_h)=h+1=2h+1-\gamma(SS_h)$. Observe that if $T_{SS}=T-SS_h$ consists of two vertices, then $T$ is a subdivided star $SS_{h+1}$ and the claim follows. If $T-SS_h$ contains only one vertex, then $u$ was a support vertex, contradicting our assumption. Furthermore, if $T_{SS}\cong P_4$, then if $u$ is adjacent to a support vertex of $P_4$, we may consider $C_h$ together with the support vertices of $P_4$ as our identifying code. If $u$ is adjacent to a leaf of $P_4$, then we may consider $C_h$ together with the leaf farthest from $u$ and the only non-leaf vertex at distance~2 from $u$. In both cases, we have $\ID(T)\le h+3=n-h-2$ since $n=2h+5$, and $\gamma(T)=h+2$. Hence, the upper bound $\ID(T)\le n-\gamma(T)$ holds. Finally, there is the case where, by induction, the tree $T_{SS}$ contains an identifying code $C_{SS}$ of size $|C_{SS}|=\ID(T_{SS})\le n-(2h+1)-\gamma(T_{SS})$. However, in this case, the set $C_{SS}\cup C_h$ is an identifying code in $T$ containing at most $n-(2h+1)-\gamma(T_{SS})+(h+1)=n-(h+\gamma(T_{SS}))\le n-\gamma(T)$ vertices.
This completes the proof.
\end{proof}

Observe that for every tree $T$ of order at least~$3$, we have $\gamma(T)\ge s(T)$ since every leaf needs to be dominated by a vertex in the dominating set. Moreover, for example, in a path, we can have $\gamma(P_n)>s(P_n)$. Therefore, the upper bound of Theorem \ref{the_domBound} is an improvement over the simple but useful upper bound $n-s(T)$ from Lemma~\ref{lem_FT_1}.

Moreover, the upper bound $n-\gamma(T)$ is tight for example for every tree in Figure \ref{fig:trees}. This was the leading motivation for us to prove this bound. If we compare the upper bound of Theorem \ref{the_main} to the upper bound $n-\gamma(T)$ of Theorem~\ref{the_domBound}, we observe that in some cases, they are equal (for example for bi-stars consisting of two equal-sized stars), but in some cases the domination bound gives a smaller value (for example for bi-stars consisting of two unequal stars). In some cases, the bound from Theorem \ref{the_main} is even smaller (for example for a chain of equal-sized stars joined with a single edge between some leaves).

\section{Characterizing extremal trees}\label{sec:small_cases}

As mentioned before, the main difficulty in proving Conjecture~\ref{conj_G Delta_UB} lies in handling the constant $c$. In this section, we focus on those extremal trees that require $c>0$, and by a careful analysis, we are able to fully characterize them. This may seem a little tedious, but the analysis is crucial for the proof of Conjecture~\ref{conj_G Delta_UB} for trees in Section~\ref{sec:main}, and also for the proof for all triangle-free graphs in the companion paper~\cite{PaperPart2}.

Since the identification number of paths is well-understood, for the rest of the section, we only consider graphs of maximum degree $\Delta \ge 3$. Towards proving the first part of Theorem~\ref{the_main}, we next look at the \emph{star} graphs. For any $\Delta \ge 3$, the complete bipartite graph $K_{1,\Delta}$ is called a \emph{$\Delta$-star}, or simply a \emph{star}. Noting that for any $\Delta$-star $S$ the set of all its leaves constitutes a minimum identifying code of $S$, it can therefore be readily verified that the following lemma is true.

\begin{lemma} \label{lem_star}
For a $\Delta$-star $S$ with $\Delta \ge 3$ of order $n ~(=\Delta+1)$, we have
\[
\ID(S) =  \left( \frac{\Delta - 1}{\Delta} \right) n  + \frac{1}{\Delta}.
\]
\end{lemma}

In particular, Lemma~\ref{lem_star} shows that $\Delta$-stars satisfy the conjectured bound with $c = \frac{1}{\Delta}$. Hence, for all stars the constant $c=\frac{1}{3}$ suffices in Conjecture~\ref{conj_G Delta_UB}. To fully establish the first part of Theorem~\ref{the_main} now, one needs to only show the veracity of the result for the rest of the trees in $\mathcal{T}_{3}$. To describe the trees in $\mathcal{T}_{3}$ and other graphs later in a more unified manner, we start by defining a particular ``join" of graphs with stars. Let $G'$ be a graph and $S$ be any star. Then, let $G' \rhd_v S$ denote the graph obtained by identifying a vertex $v$ of $G'$ with a leaf $l$ of $S$ (for example, if $S$ and $P$ are a $3$-star and $4$-path, respectively, each with a leaf $v$, then the graphs in Figures~\ref{fig:trees}(b) and \ref{fig:trees}(c) are $S \rhd_v S$ and $P \rhd_v S$, respectively). We call the $G' \rhd_v S$ the \emph{graph $G'$ appended with a star} and it is said to be obtained by \emph{appending $S$ (by its leaf $l$) onto (the vertex $v$ of) $G'$}. In the case that the vertex $v$ of $G'$ is inconsequential to the context or is (up to isomorphism) immaterial to the graph $G$, we may simply drop the suffix $v$ in the notation $G' \rhd_v S$ and denote it as $G' \rhd S$ (for example, if $P$ is a $2$-path and $S$ is a star, then $P \rhd_v S$ is (up to isomorphism) the only graph irrespective of which vertex of the $2$-path $v$ is. Another example would be if $S$ is a $\Delta$-star and we require $S \rhd_v S$ to be a graph of maximum degree $\Delta$. As a convention, we continue to call the vertices of $G' \rhd S$ by the same names as they were called in the graphs $G'$ and $S$. In other words, the graph $G' \rhd S$ is said to \emph{inherit} its vertices from $G'$ and $S$. In particular, if $G' \rhd S$ is obtained by identifying the vertex $v$ of $G'$ and a leaf $l$ of $S$, then both the names $v$ and $l$ (as and when convenient) also refer to the identified vertex in $G' \rhd S$.

\medskip

Let $G_0$ be a fixed graph, $p \ge 1$ be an integer and for each $i \in [p]$, let $S_i$ be a $\Delta_i$-star for $\Delta_i \ge 3$. Now, we may carry out the process of inductively appending stars by defining $G_i = G_{i-1} \rhd_{v_{i-1}} S_i$ for all $i \in [p]$, where $v_{i-1}$ is a vertex of $G_{i-1}$. Then the graph $G_p$ is called the \emph{graph $G_0$ appended with $p$ stars}. In the case that each $S_i$ is isomorphic to a $\Delta$-star $S$ for $\Delta \ge 3$, we call the graph $G_p$ the \emph{graph $G_0$ appended with $p$ $\Delta$-stars}. In the particular case that $G_0 = S_0$ is itself a $\Delta_0$-star for $\Delta_0 \ge 3$, we simply call $G_p$ an \emph{appended star}. Further, if $\Delta = \Delta_0 = \Delta_1 = \cdots = \Delta_p$, then we call the graph $G_p$ an \emph{appended $\Delta$-star}.

\medskip

We next furnish some general results for any identifiable graph appended with a star.

\begin{lemma}\label{lem_ID <= ID'+Delta-1}
If $G'$ is an identifiable graph and $G = G' \rhd_v S$, where $v$ is a vertex of $G'$ and $S$ is a $\Delta$-star for $\Delta \ge 3$, then $G$ is also identifiable and
\[
\ID(G) \le \ID(G') + \Delta-1.
\]
\end{lemma}

\begin{proof}
Graph $G$ is clearly identifiable. Assume now that $l_1, l_2, \ldots, l_{\Delta-1}$ are the leaves that $G$ inherits from $S$. Also, assume that $u$ is the universal vertex of $S$, that is, $u$ is the common support vertex of the leaves $l_1, l_2, \ldots , l_{\Delta-1}$ in $G$. Suppose now that $C'$ is a minimum identifying code of $G'$. Then we claim that the set $C = C' \cup \{ l_1, l_2, \ldots , l_{\Delta-1} \}$ is an identifying code of $G$. Clearly, $C$ is a dominating set of $G$, since $C'$ is a dominating set of $G'$. Next, to check that $C$ is also a separating set of $G$, it is enough to check that each of the vertices $u, l_1, l_2, \ldots , l_{\Delta-1}$ is separated from  all vertices at distance at most~$2$ in $G$ by $C$. To start with, assume $w$ to be a vertex of $G$ other than $u$. Then, one can find at least one $l_i$ to always be an separating $C$-codeword for the pair $u,w$. Moreover, each $l_i$ is itself an separating $C$-codeword for the pair $l_i,w$, where $l_i \ne w$. This establishes the claim that $C$ is an identifying code of $G$. Therefore, $\ID(G) \le |C| = |C'| + \Delta-1 = \ID(G') + \Delta - 1$. This proves the result.
\end{proof}

The next lemma shows that if $G$ is the graph obtained by starting from an identifiable ``base" graph $G_0$ and iteratively appending stars thereon, then the graphs $G$ and $G_0$ share the same constant $c$ in Conjecture~\ref{conj_G Delta_UB}.

\begin{lemma}\label{lem_ID(G) from ID(G')}
Let $c$ be a constant, $G_0$ be an identifiable graph of order $n_0$, of maximum degree $\Delta_0$ and be such that $\ID(G_0) \le ( \frac{\Delta_0-1}{\Delta_0} ) n_0 + c$. For an integer $p \ge 1$ and $i \in [p]$, let $S_i$ be a $\Delta_i$-star for $\Delta_i \ge 3$. Also, for $i \in [p]$, let $G_i = G_{i-1} \rhd_{v_{i-1}} S_i$, where $v_{i-1}$ is a vertex of $G_{i-1}$, and $G = G_p$ is a graph of order $n$ and maximum degree $\Delta$ obtained by appending $p$ stars to $G_0$. Moreover, assume that $\Delta_{\max} = \max \{\Delta_i: 0 \le i \le p \}$. Then, we have
\[
\ID(G) \le \left( \frac{\Delta_{\max}-1}{\Delta_{\max}}  \right) n + c \le \left( \frac{\Delta-1}{\Delta} \right) n + c.
\]
\end{lemma}
\begin{proof}
Let us prove the claim by induction on $q \in [0,p]$. Case $q=0$ follows from the statement itself regarding the graph $G_0$. Let us assume that the induction hypothesis holds for $q \in [0, p-1]$ and consider the case that $q = p \ge 1$. We have $G = G_{p-1} \rhd S_p$, where we let $G_{p-1}$ have order $n'$ and where $\Delta'_{\max} = \max \{ \Delta_i: 1 \le i \le p-1 \}$. Finally, by the induction hypothesis, let $C'$ be a minimum identifying code of $G_{p-1}$ of cardinality at most $( \frac{\Delta'_{\max}-1}{\Delta'_{\max}} )n' + c$. By Lemma~\ref{lem_ID <= ID'+Delta-1}, there exists an identifying code $C$ of $G$ of cardinality $|C'|+\Delta_p-1$. Therefore, we have
\begin{flalign*}
\ID(G) \le | C| = |C'|+\Delta_p-1 & \le \left( \frac{\Delta'_{\max}-1}{\Delta'_{\max}} \right) n' + c + \left( \frac{\Delta_p-1}{\Delta_p} \right) \Delta_p\\
&\le \left( \frac{\Delta_{\max}-1}{\Delta_{\max}} \right) (n-\Delta_p) + \left( \frac{\Delta_{\max}-1}{\Delta_{\max}} \right) \Delta_p + c\\
&= \left( \frac{\Delta_{\max}-1}{\Delta_{\max}} \right) n + c \\
& \le \left( \frac{\Delta-1}{\Delta} \right) n + c.
\end{flalign*}
\end{proof}

\begin{corollary} \label{cor_starbadbound}
For an integer $p \ge 1$ and $i \in [0,p]$, let $S_i$ be a $\Delta_i$-star for $\Delta_i \ge 3$. For $i \in [p]$, let $G_i = G_{i-1} \rhd_{v_{i-1}} S_i$, where $v_{i-1}$ is a vertex of $G_{i-1}$, and where $G = G_p$ is an appended star of order $n$ and of maximum degree $\Delta$. Moreover, assume that $\Delta_{\max} = \max \{\Delta_i: 0 \le i \le p \}$. Then, we have
\[
\ID(G) \le \left(  \frac{\Delta_{\max}-1}{\Delta_{\max}}  \right) n + \frac{1}{\Delta_0} \le \left(  \frac{\Delta-1}{\Delta}  \right) n + \frac{1}{3}.
\]
\end{corollary}

\begin{proof}
The result follows from taking $G_0 = S_0$ in Lemma~\ref{lem_ID(G) from ID(G')} and $c = \frac{1}{\Delta_0} \le \frac{1}{3}$ from Lemma~\ref{lem_star}.
\end{proof}

It is worth mentioning here that Lemma~\ref{lem_ID(G) from ID(G')} is central to our inductive proof arguments later, whereby, if a graph $G$ is structurally isomorphic to a ``smaller" identifiable graph $G'$ appended with a star, then by using Lemma~\ref{lem_ID(G) from ID(G')} for $p=1$ and an inductive hypothesis that the ``smaller" graph $G' = G - S$ satisfies the conjectured bound, one can show that so does the ``bigger" graph $G$. We next look at the identification numbers of the trees isomorphic to $P \rhd_v S$, where $P$ is a path, $v$ is a vertex of $P$ and $S$ is a star. In particular, we establish the identification numbers of the trees $T_2$ and $T_3$ in the collection $\mathcal{T}_{3}$.

\subsection{Paths appended with stars}\label{subsec:AppPaths}

In this subsection, we analyse the identification numbers of paths appended with stars with respect to the conjectured bound. We begin by establishing the identification numbers of the trees $T_2$ and $T_3$ in the collection $\mathcal{T}_{3}$ illustrated in Figures~\ref{fig:trees}(c) and~\ref{fig:trees}(d), respectively.

\begin{lemma} \label{lem_T2}
If $T_2$ is the tree in $\mathcal{T}_{3}$ of order~$n=7$ and of maximum degree $\Delta=3$, then
\[
\ID(T_2) = 5 = \frac{2}{3} \times 7 + \frac{1}{3} = \left( \frac{\Delta-1}{\Delta} \right) n + \frac{1}{3}.
\]
\end{lemma}
\begin{proof}
Assume that $T_2 = P \rhd_z S_1$, where $P$ is a $4$-path, $z$ is a leaf of $P$ and $S_1$ is a $3$-star (see Figure~\ref{fig:trees}(c)). Assume that $V(P) = \{w,x,y,z\}$, where $w$ and $z$ are the leaves of $P$ with support vertices $x$ and $y$, respectively. Further, assume that $V(S_1) = \{u_1, a_1, b_1, c_1\}$, where $u_1$ is the universal vertex of $S_1$ and the vertices $a_1, b_1, c_1$ are the three leaves of $S_1$. Also assume that $S_1$ is appended by its leaf $c_1$ onto the leaf $z$ of $P$.

\medskip

We first show that any identifying code of $T_2$ has order at least $5$. So, let $C$ be an identifying code of $T_2$. Then $C$ must contain at least two vertices each from the sets $W_0 = \{w,x,y\}$ and $W_1 = \{u_1,a_1,b_1\}$. Now, if either of $W_0$ or $W_1$ is a subset of $C$, then we are done. So, assume that exactly two vertices from each of $W_0$ and $W_1$ belong to $C$. Moreover, if $z \in C$, then we are done again. So, we also assume that $z \notin C$. We must have $y\in C$ as the unique separating $C$-codeword for the pair $x$ and $w$ in $G$. Moreover, for $C$ to dominate $w$, we need $C\cap\{w,x\}\ne \emptyset$. This implies that $w \in C$, or else, the pair $x$ and $y$ are not identified by $C$ in $G$ (notice that $z \notin C$ by assumption). Therefore, $C \cap W_0 = \{w,y\}$. Further, for $C$ to separate $a_1$ and $b_1$, at least one of them must be in $C$. So, without loss of generality, let us assume that $a_1 \in C$. Then, for $C$ to dominate $b_1$, we need $C \cap \{b_1,u_1\} \ne \emptyset$. Now, for $C$ to separate $y$ and $z$ in $G$, we must have $u_1\in C$ (notice that $x \notin C$). This implies that $b_1 \notin C$ (since $|C \cap W_1|=2$ by our assumption). However, this implies that $C$ does not separate $u_1$ and $a_1$ in $G$, a contradiction. Hence, $|C| \ge 5$.

\medskip
As for proving that the identifying number of $T_2$ is bounded above by $5$, since $\ID(P) = 3 = \frac{2}{3} \times 4 + \frac{1}{3}$, taking $G_0 = P$ and $c=\frac{1}{3}$ in Lemma~\ref{lem_ID(G) from ID(G')}, we have $\ID(T_2) \le \frac{2}{3} \times 7 + \frac{1}{3} = 5$.
 \end{proof}

\begin{lemma} \label{lem_T3}
If $T_3$ is the tree in $\mathcal{T}_{3}$ of order~$n=10$ and of maximum degree $\Delta=3$, then
\[
\ID(T_3) = 7 = \frac{2}{3} \times 10 + \frac{1}{3} = \left(  \frac{\Delta-1}{\Delta} \right) n + \frac{1}{3}.
\]
\end{lemma}
\begin{proof}
Assume that $T_2 = P \rhd_z S_1$ and $T_3 = T_2 \rhd_z S_2$, where $P$ is a $4$-path, $z$ is a leaf of $P$ and $S_1$ and $S_2$ are $3$-stars (see Figure~\ref{fig:trees}(d)). Assume that $V(P) = \{w,x,y,z\}$, where $w$ and $z$ are the leaves of $P$ with support vertices $x$ and $y$, respectively. Further, assume that, for $i \in [2]$, $V(S_i) = \{u_i,a_i,b_i,c_i\}$, where $u_i$ is the universal vertex of $S_i$ and $a_i, b_i, c_i$ are the three leaves of $S_i$. Assume that each $S_i$ is appended by its leaf $c_i$ onto the leaf $z$ of $P$. We first show that any identifying code of $T_3$ has order at least $7$. So, let $C$ be an identifying code of $T_3$. As in the proof of Lemma~\ref{lem_T2}, $C$ must contain at least two vertices from each of the sets $W_0 = \{w,x,y\}$ and $W_i = \{u_i,a_i,b_i\}$ for $i \in [2]$. Now, if any of $W_0$, $W_1$ and $W_2$ are subsets of $C$, then we are done. So, assume that exactly two vertices from each of $W_0$, $W_1$ and $W_2$ belong to $C$. Moreover, if $z \in C$, then we are done again. So, we also assume that $z \notin C$. Then, for both $i \in [2]$, $u_i \notin C$ (or else, exactly one of $a_i$ and $b_i$ belongs to $C$; and thus, if $a_i \in C$, without loss of generality, then $z$ is forced to be in $C$ as the only separating $C$-codeword for the pair $u_i, a_i$ in $G$, contradicting our assumption that $z \notin C$).

\medskip
On the other hand, we know that the vertex $y$ must belong to $C$ as the only separating $C$-codeword for the pair $w, x$ in $G$. Since, by our assumption, $C$ contains exactly two vertices from $W_0$, this forces exactly one of $w$ and $x$ to belong to the set $C$. Thus, if $w \notin C$, then the pair $x, y$ has no separating $C$-codeword in $G$ (since, by our assumption, $z \notin C$ either); and if $x \notin C$, then the pair $y, z$ has no separating $C$-codeword in $G$ (since, again, for both $i \in [2]$, $u_i \notin C$). Either way, we produce a contradiction. Hence, $|C| \ge 7$.

\medskip
As for proving that the identifying number of $T_3$ is bounded above by $7$, since $\ID(T_2) = 5 = \frac{2}{3} \times 7 + \frac{1}{3}$, taking $G_0 = T_2$ and $c=\frac{1}{3}$ in Lemma~\ref{lem_ID(G) from ID(G')}, we have $\ID(T_3) \le \frac{2}{3} \times 10 + \frac{1}{3} = 7$.
\end{proof}

Lemmas~\ref{lem_T2} and \ref{lem_T3} therefore establish the result in Theorem~\ref{the_main} for the trees $T_2$ and $T_3$ in the collection $\mathcal{T}_{3}$. For the rest of this subsection, we look at other paths appended with stars which are \emph{not} isomorphic to either $T_2$ or $T_3$ and show that they satisfy Conjecture~\ref{conj_G Delta_UB} with constant $c=0$.

\begin{lemma} \label{lem_(m,Delta)-broom}
Let $G = P \rhd_v S$ be a graph of order $n$, where $S$ is a $\Delta$-star with $\Delta \ge 3$, $P$ is a path and $v$ is a vertex of $P$. If any one of the following properties hold \\ [-18pt]
\begin{enumerate}
\item[{\rm (1)}] $P$ is not a $4$-path; or
\item[{\rm (2)}] $P$ is a $4$-path and $v$ is a non-leaf vertex of $P$; or
\item[{\rm (3)}] $P$ is a $4$-path, $v$ is a leaf of $P$ and $\Delta \ge 4$,
\end{enumerate}
then $\ID(G) \le ( \frac{\Delta-1}{\Delta} ) n$.
\end{lemma}
\begin{proof}
The maximum degree of $G$ is $\Delta$. Let $P$ be an $m$-path and $v$ denote the vertex of $P$ to which the star $S$ is appended. Let $l_1, l_2, \ldots , l_{\Delta-1}$ be the leaves that $G$ inherits from $S$ and let $u$ be the universal vertex of $S$, that is, $u$ is the common support vertex of the leaves $l_1, l_2, \ldots , l_{\Delta-1}$ of $G$. We note that $n = m+\Delta$. We divide the proof into cases each dealing with a type given in the statement. Furthermore, let $v$ be the vertex of $P$ adjacent to $u$ and (if they exist) the other vertices of $P$ are called,  from the closest vertex to $u$ to the farthest, $y$, $x$ and $w$.

\medskip
\noindent
\emph{Case~1. $P$ is not a $4$-path.} We further divide this case into the following subcases.

\smallskip
\noindent
\emph{Case~1.1. $m=2$.}  In this case $\{u,v\}$ is a dominating set in $P \rhd_v S$. By Theorem \ref{the_domBound}, we have
\[
\ID(P \rhd_v S)\le \Delta+2-2=\Delta < \left(  \frac{\Delta-1}{\Delta} \right) (\Delta+2) = \left(  \frac{\Delta-1}{\Delta} \right)  n.
\]
See Figure~\ref{fig_2,Delta_broom} for an illustration of an identifying code in $P \rhd_v S$.

\smallskip
\noindent
\emph{Case~1.2. $m \ge 3$ and $m \ne 4$.} Since $n=m+\Delta$, we have $(\frac{\Delta-1}{\Delta} ) n = m+\Delta-\frac{m}{\Delta}-1$. By supposition $\Delta \ge 3$. When $m\ge 3$, we have $\frac{m-1}{2}\ge \frac{m}{3}\ge \frac{m}{\Delta}$ and when $m\ge 6$, we have $\frac{m-2}{2}\ge \frac{m}{3}\ge \frac{m}{\Delta}$. Let us first consider odd $m\ge 3$. Using Lemma~\ref{lem_ID <= ID'+Delta-1} and Theorem~\ref{lem_paths_cycles}, we therefore infer that
\[
\ID(G)\le \ID(P)+\Delta-1=\frac{m+1}{2}+\Delta-1=m+\Delta-1-\frac{m-1}{2}\le m+\Delta-1-\frac{m}{\Delta} = \left(  \frac{\Delta-1}{\Delta} \right) n.
\]
We now consider even $m$. Since $m\ne 4$, we have $m\ge 6$ and hence,
\[
\ID(G)\le \ID(P)+\Delta-1=\frac{m+2}{2}+\Delta-1=m+\Delta-1-\frac{m-2}{2}\le m+\Delta-1-\frac{m}{\Delta} = \left(  \frac{\Delta-1}{\Delta} \right) n.
\]

\begin{figure}[!t]
     \centering
     \begin{subfigure}[!h]{0.3\textwidth}
         \centering
         \includegraphics[scale=0.3]{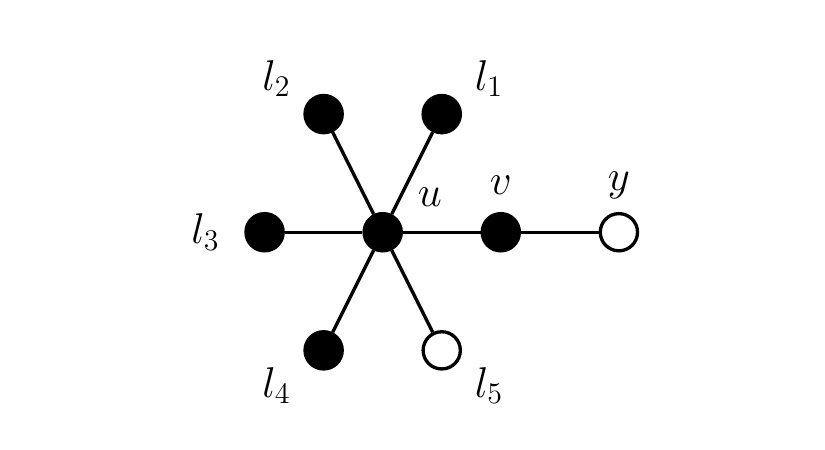}
     \caption{\small $G = P \rhd S$, where $P$ is a $2$-path and $S$ is a $6$-star.}
     \label{fig_2,Delta_broom}
     \end{subfigure}
     \hspace{2mm}
     \begin{subfigure}[!h]{0.3\textwidth}
         \centering
         \includegraphics[scale=0.3]{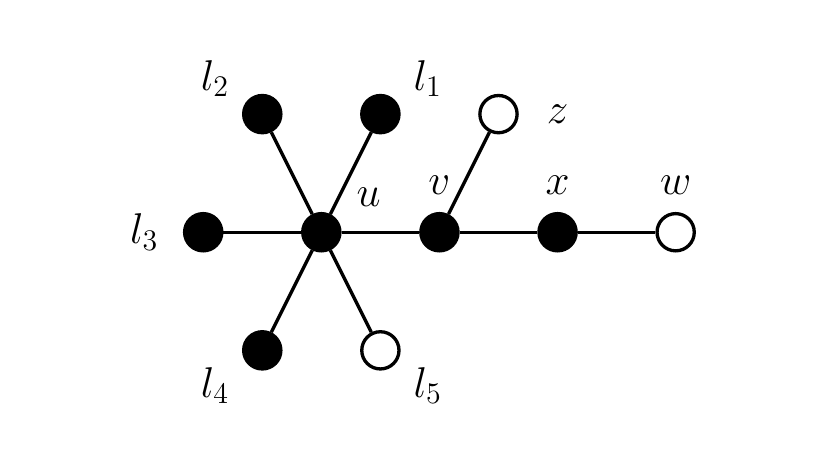}
     \caption{\small $G = P \rhd_v S$, where $P$ is a $4$-path, $v$ is a non-leaf of $P$ and $S$ is a $6$-star.}
     \label{fig_General_broom}
     \end{subfigure}
     \hspace{2mm}
     \begin{subfigure}[!h]{0.3\textwidth}
         \centering
         \includegraphics[scale=0.3]{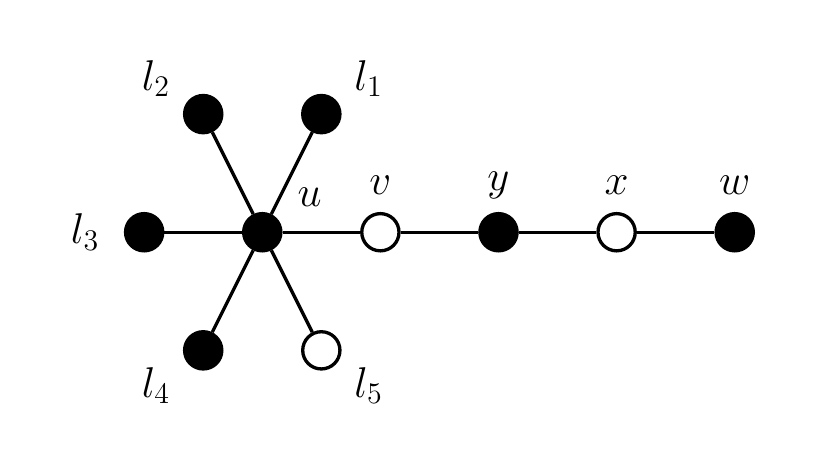}
     \caption{\small $G = P \rhd_v S$, where $P$ is a $4$-path, $v$ is a leaf of $P$ and $S$ is a $6$-star.}
     \label{fig_Perfect 4,6 broom}
     \end{subfigure}
     \caption{\small Examples for the cases in the proof of Lemma~\ref{lem_(m,Delta)-broom}. The black vertices in each graph $G$ constitute an identifying code.}
\end{figure}

\medskip
\noindent
\emph{Case~2. $P$ is a $4$-path and $S$ is appended onto one of the non-leaf vertices of $P$.}
In this case, $v$ is a non-leaf vertex of $P$. Now, let $x$ be the (other) non-leaf vertex of $P$ adjacent to $v$. The set $C = \{u, v, x, l_1, l_2, \ldots , l_{\Delta-2}\}$ is an identifying code of $G$ as illustrated in~Figure~\ref{fig_General_broom}. Since the identifying code $C$ has cardinality
\[
|C| = \Delta+1 < \left(  \frac{\Delta-1}{\Delta} \right)  (\Delta+4) = \left(  \frac{\Delta-1}{\Delta} \right) n,
\]
the desired result follows.

\medskip
\noindent
\emph{Case~3. $P$ is a $4$-path, $S$ is appended onto a leaf of $P$ and  $\Delta \ge 4$.}
In this case, $v$ is a leaf of $P$. 
The set $C = \{u, w, y, l_1, l_2, \ldots , l_{\Delta-2}\}$ is an identifying code of $G$ as illustrated in~Figure~\ref{fig_Perfect 4,6 broom}. Since the identifying code $C$ has cardinality
\[
|C| = \Delta+1 < \left(  \frac{\Delta-1}{\Delta} \right)  (\Delta+4) = \left(  \frac{\Delta-1}{\Delta} \right) n,
\]
the desired result follows. \qedhere
\end{proof}

Lemma~\ref{lem_(m,Delta)-broom} shows that all trees of the form $P \rhd S$, where $P$ is a $4$-path and $S$ is a star, but not isomorphic to $T_2$, satisfy the bound in Conjecture~\ref{conj_G Delta_UB} with $c=0$. The next corollary follows immediately from Lemma~\ref{lem_(m,Delta)-broom}.

\begin{corollary} \label{cor_(4,Delta)-broom}
If $G = P \rhd_v S$ is a graph of order $n$, where  $P$ is a path, $v$ is a vertex of $P$ and $S$ is a $\Delta$-star with $\Delta \ge 4$, then
\[
\ID(G) \le \left(  \frac{\Delta-1}{\Delta} \right) n.
\]
\end{corollary}

Our next lemma shows that all trees of the form $T_2 \rhd S$, where $S$ is a star, but not isomorphic to $T_3$ satisfy the bound in Conjecture~\ref{conj_G Delta_UB} with $c=0$.

\begin{lemma} \label{lem_T2 rhd S}
If $G = T_2 \rhd_v S$ is a graph of order~$n$ and of maximum degree $\Delta \ge 3$ such that $G \not \cong T_3$, where $v$ is a vertex of $T_2$ and $S$ is a $\Delta_S$-star with $\Delta_S \ge 3$, then
\[
\ID(G) \le \left(  \frac{\Delta-1}{\Delta} \right) n.
\]
\end{lemma}

\begin{figure}[!t]
     \centering
     \begin{subfigure}[!h]{0.48\textwidth}
         \centering
         \includegraphics[scale=0.3]{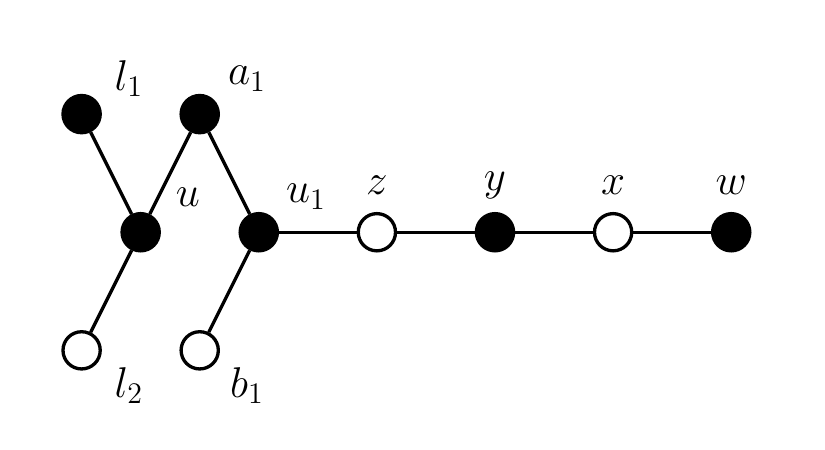}
     \caption{\small $G = T_2 \rhd_{a_1} S$, where $S$ is a $3$-star and $a_1$ is a leaf of $T_2$.}
     \label{fig_W1 rhd S_2}
     \end{subfigure}
     \hspace{2mm}
     \begin{subfigure}[!h]{0.48\textwidth}
         \centering
         \includegraphics[scale=0.3]{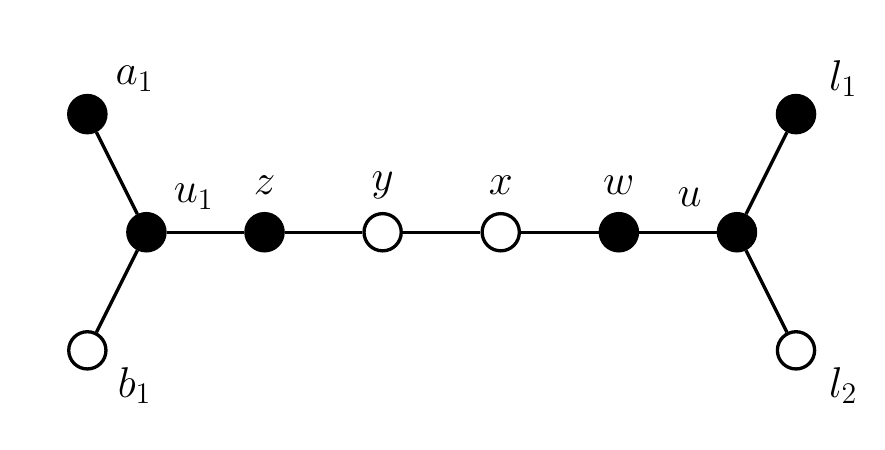}
     \caption{\small $G = T_2 \rhd_w S$, where $S$ is a $3$-star and $w$ is a leaf of $T_2$.}
     \label{fig_W1 rhd S_3}
     \end{subfigure}
     \caption{\small Figures illustrating the cases in Lemma~\ref{lem_T2 rhd S}. The black vertices constitute an identifying code.}
     \label{fig_W1 rhd S}
\end{figure}

\begin{proof}
Clearly, $\Delta_S \le \Delta$. We consider the following cases.

\medskip
\noindent
\emph{Case~1. $\Delta \ge 4$.} Observe in this case that either $\Delta_S=\Delta$ or $\Delta=4$ and $\Delta_S=3$. By Lemma~\ref{lem_T2}, we have $\ID(T_2)=5$. Therefore, by Lemma~\ref{lem_ID <= ID'+Delta-1} we infer that
\[
\ID(G)\le \ID(T_2)+\Delta_S-1=\Delta_S+4 < \left(  \frac{\Delta-1}{\Delta} \right) (\Delta_S+7) = \left(  \frac{\Delta-1}{\Delta} \right) n,
\]
yielding the desired result in this case.

\medskip
Thus, we assume from now on that $\Delta=\Delta_S=3$. Hence, the graph $G$ has order $n=10$ and it suffices for us to show that there exists an identifying code $C$ of $G$ of cardinality $|C| = 6<\frac{2}{3} \times 10 = ( \frac{\Delta-1}{\Delta} ) n$. Let $T_2 = P \rhd_z S_1$, where $P$ is a $4$-path, $z$ is a leaf of $P$ and $S_1$ is a $\Delta$-star for $\Delta \ge 3$. Assume further that \\ [-18pt]
\begin{enumerate}
\item[{\rm (1)}] $V(P)=\{w,x,y,z\}$, where $w$ and $z$ are leaves of $P$ with support vertices $x$ and $y$, respectively. \1
\item[{\rm (2)}] $V(S_1)=\{u_1, a_1, b_1, c_1\}$, where $u_1$ is the universal vertex of $S_1$, the vertices $a_1, b_1, c_1$ are the three leaves of $S_1$ and that $S_1$ is appended by its leaf $c_1$ onto the leaf $z$ of $P$. \1
\item[{\rm (3)}] $V(S)=\{u, l_1, l_2, l_3\}$, where $u$ is the universal vertex of $S$, the vertices $l_1, l_2, l_3$ are the leaves of $S$ and $l_3$ is the leaf of $S$ by which the latter is appended onto the vertex $v$ of $T_2$.
\end{enumerate}

Since $\Delta=\Delta_S=3$, $S$ is not appended onto the universal vertex $u_1$ of $S_1$. We consider next two further cases.

\medskip
\noindent
\emph{Case~2. $S$ is appended onto either of the leaves $a_1$ and $b_1$ of $S_1$.} Without loss of generality, let us assume that $S$ is appended onto the leaf $a_1$ of $T_2$. The set $C=\{y,w,u_1,u,a_1,l_1\}$ is an identifying code of $G$ of cardinality~$6$ as illustrated in~Figure~\ref{fig_W1 rhd S_2}, where the vertices in the identifying code are marked with black, implying that $\ID(G) \le 6 < \frac{2}{3} \times 10 = ( \frac{\Delta-1}{\Delta} ) n$, yielding the desired result.

\medskip
\noindent
\emph{Case~3. $S$ is appended onto any vertex of $T_2$ other than $u_1$, $a_1$ and $b_1$.}
In this case, since $G \not \cong T_3$, the star $S$ cannot be appended onto the vertex $z$ of $T_2$ (that is, $v \ne z$). Thus, $S$ is appended onto either of the vertices $w$, $x$ and $y$ of $T_2$. We further divide this case into the following subcases.

\smallskip
\noindent
\emph{Case~3.1. $S$ is appended onto either of the non-leaf vertices $x$ or $y$ of $T_2$.}
In this case, the graph $G$ can also be expressed as $G = G' \rhd_z S_1$, where $G' = P \rhd_v S$.
Since the maximum degree of $G'$ is $\Delta_S=3$, by Lemma~\ref{lem_(m,Delta)-broom}(2), we have
\[
\ID(G') \le \frac{2}{3} |V(G')| = \frac{2}{3} (n-3).
\]
The bound for $G$ therefore follows with Lemma~\ref{lem_ID <= ID'+Delta-1}.

\smallskip
\noindent
\emph{Case~3.2. $S$ is appended onto the leaf $w$ of $T_2$.}
In this case, the set $C = \{a_1, u_1, z,$ $w, u, l_1\}$ is an identifying code of $G$ of cardinality~$6$ as illustrated in~Figure~\ref{fig_W1 rhd S_3} where the vertices in the identifying code are marked with black, implying that $\ID(G) \le 6 < ( \frac{\Delta-1}{\Delta} ) n$, as desired. \qedhere
\end{proof}

Finally, the next lemma shows that all trees of the form $T_3 \rhd S$, where $S$ is a star, also satisfy the bound in Conjecture~\ref{conj_G Delta_UB} with $c=0$.

\begin{lemma} \label{lem_T3 rhd S}
If $G = T_3 \rhd_v S$ is a graph of order~$n$ and of maximum degree $\Delta \ge 3$, where $v$ is a vertex of $T_3$ and $S$ is a $\Delta_S$-star with $\Delta_S \ge 3$, then
\[
\ID(G) \le \left(  \frac{\Delta-1}{\Delta} \right) n.
\]
\end{lemma}

\begin{proof}
Let $P$ be a $4$-path and let $S_1$ and $S_2$ be $3$-stars such that $V(P)=\{w,x,y,z\}$, where $w$ and $z$ are leaves of $P$ with support vertices $x$ and $y$, respectively. Moreover, for $i \in [2]$, we have $V(S_i)=\{u_i, a_i, b_i, c_i\}$, where $u_i$ is the universal vertex of $S_i$, and the vertices $a_i, b_i, c_i$ are the three leaves of $S_i$. Also assume that each $S_i$ is appended by its leaf $c_i$ onto the leaf $z$ of $P$ to form the tree $T_3$. Now, let $T_2 \cong P \rhd_{z} S_1$ and we have $T_3\cong (P \rhd_z S_1) \rhd_z S_2$. Let us assume first that $S$ is also appended onto the vertex $z$ of $P$ (that is, $v = z$). In this case, $\deg_G(z)=4$. If $\Delta_S=3$, then $n=13$ and by Lemma~\ref{lem_ID <= ID'+Delta-1}, we have
\[
\ID(G)\le \ID(T_3)+2=9 < \frac{3}{4}\times 13 = \left(  \frac{\Delta-1}{\Delta} \right) n.
\]
This implies that $S$ is a $\Delta$-star with $\Delta \ge 4$. Let $G' = P \rhd_z S$ have order $n' ~(=n-6)$ and  maximum degree $\Delta' \le \Delta$. By Corollary~\ref{cor_(4,Delta)-broom}, we have
\[
\ID(G') \le \left( \frac{\Delta'-1}{\Delta'} \right) n' \le \left( \frac{\Delta-1}{\Delta} \right)  (n-6).
\]
Therefore, by Lemma~\ref{lem_ID(G) from ID(G')}, we have the result.

So, let us assume that $G'' = T_2 \rhd_v S$ has order $n'' ~(= n-3)$ and maximum degree $\Delta'' \le \Delta$, where $S$ is appended onto a vertex of $T_2$ other than $z$ (that is, $v \ne z$). In that case, the graph $G'' \not \cong T_3$ and, therefore, by Lemma~\ref{lem_T2 rhd S}, we have
\[
\ID(G'') \le \left(   \frac{\Delta''-1}{\Delta''} \right) n'' \le \left(   \frac{\Delta-1}{\Delta} \right) (n-3).
\]
Therefore, we have the result by Lemma~\ref{lem_ID(G) from ID(G')}.
\end{proof}

To conclude the subsection, we close with the following remark.

\begin{remark}
If $G$ is a graph isomorphic to a path appended with stars, then $G$ satisfies Conjecture~\ref{conj_G Delta_UB} with $c=\frac{1}{3}$ in the case of $T_2$ and $T_3$, and with $c=0$ otherwise.
\end{remark}

\subsection{Appended stars}\label{Subsec:stars}

The trees in the collection $\mathcal{T}_{3}$ other than $T_2$ and $T_3$ are precisely the appended $3$-stars of maximum degree~$3$ and of diameter at most~$6$. In this subsection, we show that all appended stars satisfy the conjectured bound with $c=0$, except for those in $\mathcal{T}_{3}$, which satisfy the bound with $c=\frac{1}{3}$. To start with, we first look at appended stars of maximum degree at least~$4$.

\begin{lemma} \label{lem_4-stars}
If $G$ is an appended star of order~$n$ and of maximum degree $\Delta \ge 4$, then
\[
\ID(G) \le \left(  \frac{\Delta-1}{\Delta} \right) n.
\]
\end{lemma}
\begin{proof}
Let $p$ be an integer and, for $i \in [0,p]$, let $S_i$ be a $\Delta_i$-star with $\Delta_i \ge 3$. Also, for $i \in [p]$, let $G_i = G_{i-1} \rhd_{v_{i-1}} S_i$, where $v_{i-1}$ is a vertex of $G_{i-1}$, and such that $G = G_p$ is the appended star. Let $V(S_i)=\{u_i,l_1^i,\dots,l_{\Delta_i}^i\}$, where $u_i$ is the universal vertex and each $l^i_j$, for $j \in [\Delta_i]$, is a leaf of $S_i$. For $h \in [p]$, let $G_h$ have order $n'$ and maximum degree $\Delta'$.  Observe that if $\ID(G_h)\le ( \frac{\Delta'-1}{\Delta'} ) n'$, then $\ID(G)\le (  \frac{\Delta-1}{\Delta} ) n$ by Lemma~\ref{lem_ID(G) from ID(G')}.
We next look at the following two cases according to whether the $\Delta_i$'s are all equal or not.

\medskip
\noindent
\emph{Case~1. There exists $1 \le h \le p$ such that $\Delta_0 = \Delta_1=\cdots=\Delta_{h-1}\ne \Delta_{h}$.}
Let $G_{h-1}=G''$ have order $n''$ and maximum degree $\Delta''$. Moreover, let $C''$ be a minimum identifying code of $G''$.

\smallskip
\noindent
\emph{Case~1.1. $\Delta_h=\Delta_{h-1}+q$ for some $q>0$.} Observe that
\[
n'' = h\Delta_0+1 \mbox{ and } n'=n''+\Delta_h=(h+1)\Delta_h-hq+1.
\]
By Corollary~\ref{cor_starbadbound}, we have
\[
|C''|\le \left( \frac{\Delta_0-1}{\Delta_0} \right) n'' + \frac{1}{\Delta_0}= h(\Delta_0-1)+1=h(\Delta_h-q-1)+1.
\]
Moreover, by Lemma~\ref{lem_ID <= ID'+Delta-1}, we have
\begin{align*}
\ID(G_h) \le |C''|+\Delta_h-1 &\le h(\Delta_h-q-1)+1+(\Delta_h-1)\\
    & = \left(  \frac{\Delta_h-1}{\Delta_h} \right) h\Delta_h- hq + \Delta_h\\
    &= \left( \frac{\Delta_h-1}{\Delta_h} \right) (n' - \Delta_h+ hq -1)- hq +\Delta_h\\
    &= \left( \frac{\Delta_h-1}{\Delta_h} \right) n' - (hq - 1) \frac{1}{\Delta_h} \\
    & \le \left( \frac{\Delta'-1}{\Delta'} \right) n'.
\end{align*}

\smallskip
\noindent
\emph{Case~1.2. $\Delta_h=\Delta_{h-1}-q$ for some $q>0$.}
Observe that
\[
n''=h\Delta_0+1 \mbox{ and } n'=n''+\Delta_h=(h+1)\Delta_0-q+1.
\]
By Corollary~\ref{cor_starbadbound}, we have
\[
|C''|\le \left( \frac{\Delta_0-1}{\Delta_0} \right) n''  + \frac{1}{\Delta_0}= h(\Delta_0-1)+1.
\]
Moreover, by Lemma~\ref{lem_ID <= ID'+Delta-1}, we have
\begin{align*}
\ID(G_h) \le |C''|+\Delta_h-1 &\le h(\Delta_0-1) + \Delta_0-q\\
    &= \left(  \frac{\Delta_0-1}{\Delta_0} \right) h \Delta_0 + \Delta_0 -q\\
    &= \left(  \frac{\Delta_0-1}{\Delta_0} \right) (n' - \Delta_0 +q-1) + \Delta_0 -q\\
    &= \left(  \frac{\Delta_0-1}{\Delta_0}\right) n' - \frac{q-1}{\Delta_0}\\
    & \le \left(  \frac{\Delta_0-1}{\Delta_0}\right) n' \le \left( \frac{\Delta'-1}{\Delta'} \right) n'.
    %
    \end{align*}

In both the above subcases therefore, the result follows from Lemma~\ref{lem_ID(G) from ID(G')}. This completes the proof of Case~1.

\medskip
\noindent
\emph{Case~2.  $\Delta_0 = \Delta_1 = \cdots = \Delta_p$.} In this case, we have $n=(p+1)\Delta_0+1$ and $n'=(h+1)\Delta_0+1$. Let us now divide our analysis into the following subcases.

\smallskip
\noindent
\emph{Case~2.1. For some $i \in [p]$, the vertex $v_{i-1} \in \{ u_0, u_1, \ldots , u_{i-1} \}$.}
This subcase implies that the star $S_i$ is appended onto a universal vertex of any of the stars $S_0, S_1, \ldots , S_{i-1}$. Then by a possible renaming of the stars, let us assume that $v_0 = u_0$ is the universal vertex of $S_0$. We then take $h=1$ and look at the graph $G_1 = S_0 \rhd_{u_0} S_1$. Therefore, we have $\Delta'=\Delta_0+1 \ge 4$. The set $C'$ consisting of all the leaves of $G_1$ is an identifying code of $G_1$. Therefore, we have

\[
\ID(G_1) \le |C'| = 2\Delta_0 - 1 = \left(  \frac{\Delta_0}{\Delta_0 + 1} \right) (2\Delta_0 + 1) - \frac{1}{\Delta_0 + 1} \le \left(  \frac{\Delta' - 1}{\Delta'} \right) n'.
\]
Replacing $G_0$ by $G_1$ in Lemma~\ref{lem_ID(G) from ID(G')}, yields the desired result.

\smallskip
\noindent
\emph{Case~2.2. For each $1 \le i \le p$, the star $S_i$ is appended onto a leaf of $G_{i-1}$.}
This subcase can be further divided according to whether $\Delta_0 \ge 4$ or $\Delta_0 = 3$.

\smallskip
\noindent
\emph{Case~2.2.1. $\Delta_0 \ge 4$.}
We again take $h=1$ and consider the graph $G_1$ of order $n'=2\Delta_0+1$ with $\Delta' = \Delta_0 \ge 4$. In this case, $G_1$ has $2\Delta'-2$ leaves, a single degree~$2$ vertex, namely $l_1^2$, and two degree $\Delta'$ vertices, namely $u_1$ and $u_2$. The set $C'=\{u_1,u_2,l_2^1,\dots ,l_{\Delta'-1}^1,l_2^2,\dots ,l_{\Delta'-1}^2\}$ is an identifying code of $G_1$. Moreover, we have
\[
|C'|=2(\Delta'-1) = \left(\frac{\Delta'-1}{\Delta'} \right) n' - \frac{\Delta'-1}{\Delta'} < \left(\frac{\Delta'-1}{\Delta'} \right) n'.
\]
Therefore, again replacing $G_0$ by $G_1$ in Lemma~\ref{lem_ID(G) from ID(G')}, the result follows.

\smallskip
\noindent
\emph{Case~2.2.2. $\Delta_0=3$.}
Since $\Delta \ge 4$, this implies, again by a possible renaming of the stars, that the graph $G_3$ is obtained by appending the stars $S_1$, $S_2$ and $S_3$ onto a single leaf $l^0_1$, say, of $S_0$. The set $C' = \{ l^0_1, l^0_2, l^1_1, l^2_1, l^3_1, u_0, u_1, u_2, u_3 \}$ is an identifying code of $G_3$. This implies that
\[
|C'| = 9 < \frac{3}{4} \times 13 \le \left( \frac{\Delta'-1}{\Delta'} \right) n'.
\]
Hence, again replacing $G_0$ by $G_3$ in Lemma~\ref{lem_ID(G) from ID(G')}, the result follows. \qedhere
\end{proof}

The preceding lemma shows that if we require $c>0$ in the conjectured bound for an appended star, then the appended star must be subcubic. For the rest of this subsection, therefore, we investigate the veracity of Conjecture~\ref{conj_G Delta_UB} for appended $3$-stars of maximum degree~$3$. Before that, one can verify the following.

\begin{remark} \label{rem_star append diam}
The diameter of a subcubic appended $3$-star is even.
\end{remark}

We now look at the appended $3$-stars in $\mathcal{T}_{3}$ in the following proposition, which proves that the trees in $\mathcal{T}_{3}$ satisfy the conjectured bound with $c = \frac{1}{3}$.

\begin{proposition}\label{prop_appended stars <= 6}
If $G$ is an appended $3$-star of order $n$, of maximum degree $\Delta = 3$ and of diameter at most~$6$, then
\[
\ID(G) = \frac{2}{3}n+\frac{1}{3}.
\]
\end{proposition}
\begin{proof}
For $\diam(G) = 2$, that is, $G$ being a $3$-star itself, the result follows from Lemma~\ref{lem_star}. We  therefore consider $\diam(G) > 2$, that is, by the parity of $\diam(G)$ by Remark~\ref{rem_star append diam}, only the cases where $\diam(G)$ is either $4$ or $6$. We first show that $\ID(G) \ge \frac{2}{3} n + \frac{1}{3}$. Let $p$ be an integer and for each $i \in [0,p]$, let $S_i$ be a $3$-star. Also, for $i \in [p]$, let $G_i = G_{i-1} \rhd_{v_{i-1}} S_i$, where $v_{i-1}$ is a vertex of $G_{i-1}$, and where $G = G_p$ is the appended $3$-star.

Suppose, to the contrary, that $G$ (with diameter either $4$ or $6$) is a counterexample of minimum order such that $\ID(G) < \frac{2}{3} n + \frac{1}{3}$. Since $\diam(G)$ is either $4$ or $6$, there exists a $3$-star $S_0$, say, such that $G$ is obtained by appending other $3$-stars to the leaves of $S_0$ only (notice that, for $G$ to be subcubic, no $3$-star $S_j$ is appended onto the universal vertex of another $3$-star $S_i$). For each $i \in [0,p]$, let $V(S_i)=\{u_i, a_i, b_i, c_i\}$, where $u_i$ is the universal vertex and $a_i, b_i, c_i$ are the leaves of $S_i$. Moreover, assume that each $S_i$ is appended by its leaf $c_i$ onto the graph $G_{i-1}$. Then, for all $i \in [p]$, $a_i$ and $b_i$ are leaves of $G$ adjacent to the common support vertex $u_i$ (see Figure~\ref{fig_diam <= 6} for the tree $T_7 \in \mathcal{T}_{3}$ as a sample appended $3$-star of diameter at most~$6$; also refer to Figure~\ref{fig:trees}(a) and Figures~\ref{fig:trees}(d)-\ref{fig:trees}(k) for other such examples). Now, assume $C$ to be a minimum identifying code of $G$. Then, we claim the following.

\begin{figure}[t!]
     \centering
         \includegraphics[scale=0.4]{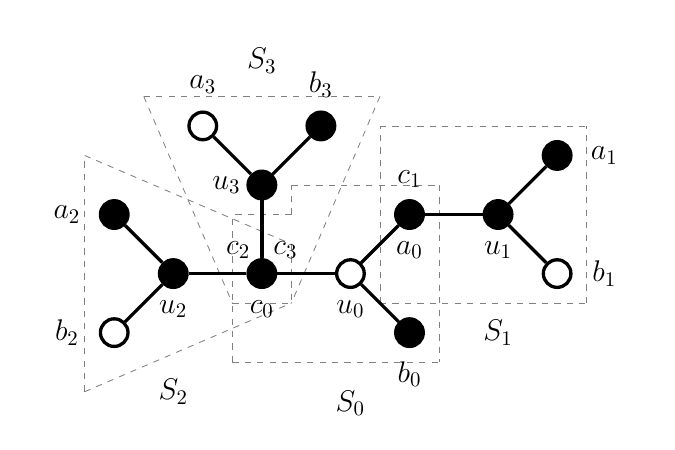}
         \caption{\small Tree $T_7 \in \mathcal{T}_{3}$: $\diam(T_7) = 6$ and $\ID(T_7) = 9$. The set of black vertices constitutes an identifying code of $T_7$.}
         \label{fig_diam <= 6}
\end{figure}

\begin{claim}\label{clm:univ-3star}
For each $i \in [p]$, the universal vertex $u_i$ of the $3$-star $S_i$ belongs to $C$.
\end{claim}
\begin{proofofclaim}
Suppose, to the contrary, that $u_i \notin C$ for some $i \in [p]$. For the vertices $a_i$ and $b_i$, both of which have degree~$1$ in $G$, to be dominated by $C$, we must have $a_i, b_i \in C$. Now, observe that, by a possible renaming of the stars $S_1, S_2, \ldots , S_p$, we can assume that $u_p \notin C$ and that $a_p, b_p \in C$. Let $C' = C \setminus \{a_p,b_p \}$. Then, $C'$ is an identifying code of $G_{p-1}$. Noting that $n = 3p+4$ and $|C| < \frac{2}{3} n + \frac{1}{3}$, we have
\begin{flalign*}
&|C'| = |C| - 2 < \frac{2}{3} n - \frac{5}{3} = \frac{2}{3} (3p+4) - \frac{5}{3} = 2p+1\\
\implies &|C'| \le 2p < 2p+\frac{2}{3} = \frac{2}{3} (3p+1) = \frac{2}{3} (n-3).
\end{flalign*}
However, since $G_{p-1}$ has order $n-3$, this contradicts the minimality of the order of $G$. Thus, we must have $u_i \in C$ for all $i \in [p]$. \end{proofofclaim}

Claim~\ref{clm:univ-3star} implies that, for each $i \in [p]$, at least two vertices from each set $\{a_i,b_i,c_i\}$, in addition to $u_i$, must belong to $C$ (otherwise the vertices of $S_i$ are not separated). However, if both $a_i, b_i \in C$, then we can discard $b_i$ from $C$ and include $c_i$ (if not included already) in the identifying code $C$ of $G$. In this way, $C$ still remains an identifying code of $G$ and of order at most the same as before. In particular then, for each $i \in [p]$, we can assume that $\{ u_i, a_i, c_i \} \subset C$. Next, we prove the following claim.

\begin{claim}
\label{claim:lower-bound}
$|C| \ge 2p+3$.
\end{claim}
\begin{proofofclaim}
If $\diam(G) = 4$, then in this case, $p \le 2$ and onto \emph{exactly one} leaf $c_0$, say, of $S_0$ at least one other $3$-star $S_1$, say, is appended by its leaf $c_1$. However then at least two vertices from $\{u_0, a_0, b_0\}$ must be in $C$. By the above considerations following the proof of Claim~\ref{clm:univ-3star}, we have at least $2p+1$ further vertices in $C$, and thus $|C| \ge 2p+3$. Hence, we may assume that $\diam(G) = 6$, for otherwise the desired lower bound in the claim follows. In this case, $p \le 6$ and to \emph{at least two} leaves $a_0$ and $c_0$ of $S_0$, other $3$-stars are appended.
Then again, by our previous observation, $a_0$ and $c_0$ are included in the identifying code $C$, since $a_0=c_i$ for some $i\ne 0$ and $c_0=c_j$ for some $j \ne 0$. Moreover, for $b_0$ to be dominated by $C$, at least one of $\{u_0,b_0\}$ must be in $C$. Thus, we again have $|C| \ge 2p+3$. This proves the claim.
\end{proofofclaim}

Finally, again noting that $n = 3p+4$, we have by Claim~\ref{claim:lower-bound} that
\[
\ID(G) = |C| \ge 2p+3 = \frac{2}{3}(3p+3)+1 = \frac{2}{3} n+ \frac{1}{3},
\]
contradicting our initial supposition that $\ID(G) < \frac{2}{3} n + \frac{1}{3}$. Hence, $\ID(G) \ge \frac{2}{3} n + \frac{1}{3}$. On the other hand, by Corollary~\ref{cor_starbadbound}, we have $\ID(G) \le \frac{2}{3}n + \frac{1}{3}$. Consequently, $\ID(G) = \frac{2}{3} n + \frac{1}{3}$.
\end{proof}

\begin{figure}[t!]
     \centering
         \includegraphics[scale=0.4]{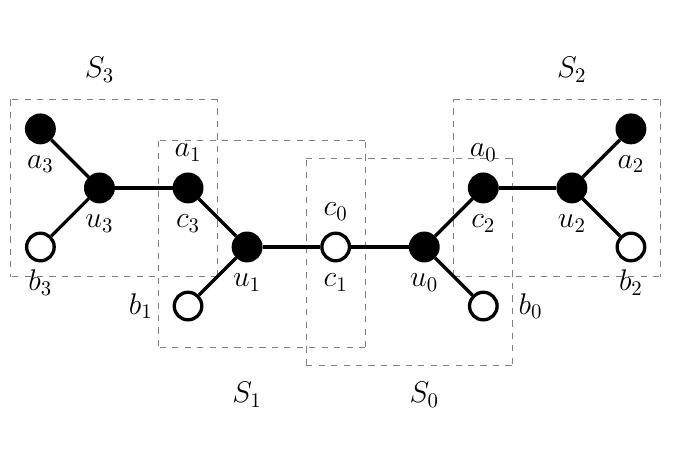}
         \caption{\small Graph $Z$: $\ID(Z) = 8$. The set of black vertices constitutes an identifying code of $Z$.}
         \label{fig_Z}
\end{figure}

We now turn our attention to the appended $3$-stars of diameter strictly larger than~$6$, that is, at least~$8$ (by the parity of the diameter in Remark~\ref{rem_star append diam}). To that end, we first show in the next lemma the existence of a special appended $3$-star of diameter~$8$.

\begin{lemma}\label{lem_uniqueZ}
Let each of $S_0, S_1, S_2, S_3$ be a $3$-star. Let $G_0 = S_0$ and, for $i \in [3]$, let $G_i = G_{i-1} \rhd_{v_{i-1}} S_i$, where $v_{i-1}$ is a vertex of $G_{i-1}$. Finally, let $Z = G_3$ be a subcubic appended $3$-star of diameter~$8$. Then, up to isomorphism, the graph $Z$ is unique.
\end{lemma}
\begin{proof}
Since $Z = G_2 \rhd_{v_2} S_3$ and $\diam(Z) = 8$, it follows that $\diam(G_2) = 6$. Next, to obtain $G_2$, assume, without loss of generality, that $S_1$ and $S_2$ are appended onto two distinct leaves of $S_0$. Then, the longest path in $G_2$ (with path-length of $\diam(G_2)$) has one leaf each of $S_1$ and $S_2$ as its endpoints. Thus, $S_3$ is appended onto a leaf of $G_2$ that is also a leaf of either $S_1$ or $S_2$ (see Figure~\ref{fig_Z}). This completes the proof.
\end{proof}

In the following lemma we establish the identification number of the graph $Z$ constructed in the statement of Lemma~\ref{lem_uniqueZ}.

\begin{lemma} \label{lem_Z}
If $Z$ is the unique subcubic appended $3$-star of diameter~$8$ as defined in the statement of Lemma~\ref{lem_uniqueZ}, then $\ID(Z)=8 < \frac{2}{3} |V(Z)|$.
\end{lemma}
\begin{proof}
For $0 \le i \le 3$, let $V(S_i) = \{u_i,a_i,b_i,c_i\}$ such that $u_i$ is the universal vertex and $a_i,b_i,c_i$ are the three leaves of $S_i$. To obtain $Z$ (up to isomorphism), let us assume that the $3$-stars are appended in such a manner that the following pairs of vertices are identified: $(c_0, c_1)$, $(a_0,c_2)$ and $(a_1,c_3)$ (see Figure~\ref{fig_Z}). 
To first show that $\ID(Z) \le 8$, it is enough to check that the set $\{u_0,u_1,u_2,u_3,a_0,a_1,a_2,a_3\}$ is an identifying code of $Z$ (again, see Figure~\ref{fig_Z} for the identifying code represented by the black vertices). We next show the reverse that $\ID(Z) \ge 8$. Let $C$ be a minimum-size identifying code in $Z$.

\setcounter{claim}{0}

\begin{claim}
If $c_0\not\in C$, then $|C|\ge 8$.
\end{claim}
\begin{proofofclaim}
Let $G'=G-\{c_0\}$. Observe that $G'$ consists of two identical components of order~$6$. We show that each of them has at least four vertices in the identifying code. First of all, $a_0\in C$ is forced as the only vertex which separates $u_0$ and $b_0$.
Moreover, $|\{u_0,b_0\}\cap C|\ge 1$ since $C$ dominates $b_0$. Finally, we have $|\{a_2,u_2,b_2\}\cap C|\ge 2$ to dominate and separate $a_2$ and $b_2$. Thus, $|C|\ge 8$.
\end{proofofclaim}

\begin{claim}
If $c_0\in C$, then $|C|\ge 9$.
\end{claim}
\begin{proofofclaim}
We again show that $|C\cap \{u_0,b_0,a_0,u_2,a_2,b_2\}|\ge 4$ and the claim follows from symmetry. If $u_2\not\in C$, then in this case, we need $a_2,b_2\in C$ to dominate them and we also need two vertices from $\{a_0,u_0,b_0\}$. Indeed, a single vertex can dominate all of these three vertices only if it is $u_0$. However, $u_0$ does not separate $a_0$ and $b_0$, implying that $|C| \ge 9$. On the other hand, if $u_2\in C$, then in this case, we have $a_2$ or $b_2$ in $C$. We may assume without loss of generality that $a_2\in C$. To separate $u_2$ and $a_2$, we also need $a_0$ or $b_2$ in $C$. Thus, $|C\cap\{a_2,b_2,u_2,a_0\}|\ge 3$. Finally, we need $u_0$ or $b_0$ in $C$ to dominate $b_0$, once again implying that $|C|\ge 9$.
\end{proofofclaim}

The lemma follows from these two claims.
\end{proof}

The above lemma shows that the graph $Z$ satisfies Conjecture~\ref{conj_G Delta_UB} with $c = 0$. We end this section with the following lemma which shows that all appended $3$-stars of maximum degree~$3$ other than those in the collection $\mathcal{T}_{3}$ (that is, with diameter at least $8$) satisfy the conjectured bound with $c=0$.

\begin{lemma} \label{lem_3-star diam 8}
If $G$ is an appended star of order $n$, of maximum degree $\Delta=3$ and of diameter at least $8$, then
\[
\ID(G) \le \frac{2}{3} n  = \left( \frac{\Delta-1}{\Delta} \right) n.
\]
\end{lemma}

\begin{proof}
Since $\diam(G) \ge 8$, the graph $G$ must contain the graph $Z$ (of Lemma~\ref{lem_Z}) as an induced subgraph. Taking $G_0 = Z$ in Lemma~\ref{lem_ID(G) from ID(G')}, the result follows.
\end{proof}

We summarize the current subsection with the following remark.

\begin{remark}\label{rem_summary stars}
If $G$ is an appended star, then $G$ satisfies Conjecture~\ref{conj_G Delta_UB} with $c = \frac{1}{3}$ if $G$ is isomorphic to a tree in $\mathcal{T}_{3}$ and with $c = 0$ otherwise.
\end{remark}

We also state the following proposition, which will be useful in our proofs.

\begin{proposition}\label{prop_Treecodes}
If $T$ is a tree in $\mathcal{T}_{3}$, then the following properties hold.  \\ [-18pt]
\begin{enumerate}
\item[{\rm (i)}]  If $T\ne T_2$, then $T$ has an optimal identifying code $C(T)$ containing all vertices of degree at most~$2$.
\item[{\rm (ii)}]  If $T\notin\{T_2,T_3\}$, then $C(T)$ can be chosen as an independent set which contains every vertex of degree at most~$2$. When we delete any code vertex $v$ from $T$, set $C(T)\setminus\{v\}$ forms an optimal identifying code of the forest $T-v$.
\end{enumerate}
 \end{proposition}
 \begin{proof}
  The  constructions in the claim can be confirmed by the identifying codes provided in Figure~\ref{fig:trees}. Their optimality follows from Lemmas~\ref{lem_star}, \ref{lem_T2}, \ref{lem_T3} and Proposition~\ref{prop_appended stars <= 6}.
 \end{proof}

\section{Proof of Theorem~\ref{the_main}}\label{sec:main}

We are now ready to prove Theorem~\ref{the_main}. The proof first utilizes Lemmas~\ref{lem_FT_1} and~\ref{lem_FT_2}, using which, we can show that one needs to consider only the trees $G$ of the form $G' \rhd_v S$,  whereby $G'$ is necessarily a tree as well. Induction plays a central role in proving Theorem~\ref{the_main}.

\begin{proof}[Proof of Theorem~\ref{the_main}]
The proof is by induction on the order $n$ of the tree $G$. Since we have $\Delta \ge 3$, this implies that $n\ge 4$. However, $n=4$ implies that $G$ is a $3$-star and thus is isomorphic to a graph in $\mathcal{T}_3$. Therefore, we take the base case of the induction hypothesis to be when $n = 5$. In the base case now, for $G$ to be a tree and non-isomorphic to a star, we must have $G \cong P \rhd S$, where $P$ is a $2$-path and $S$ is a $3$-star. Therefore, by Lemma~\ref{lem_(m,Delta)-broom}(1), the result is true in the base case. Let us assume that the induction hypothesis is true for all trees $G''$ on $n''$ vertices such that $5 \le n'' < n$, not isomorphic to a tree in $\mathcal{T}_{\Delta''}$, of maximum degree $\Delta'' \ge 3$. Let $\ell$ and $s$ be the number of leaves and support vertices, respectively, in $G$. First of all, if $s \ge \frac{n}{\Delta}$, then, by Lemma~\ref{lem_FT_1}, we have
\[
\ID(G) \le n-s \le n- \frac{1}{\Delta} n = \left( \frac{\Delta-1}{\Delta} \right) n
\]
and, hence, we are done. Moreover, if $\ell \le ( \frac{\Delta - 2}{\Delta} ) n$, then, by Lemma~\ref{lem_FT_2}, we have
\[
\ID(G) \le \frac{n+\ell}{2} \le \Big(1+\frac{\Delta-2}{\Delta}\Big) \frac{n}{2} = \left( \frac{\Delta-1}{\Delta} \right) n
\]
and we are done in this case too. We therefore assume that both $s < \frac{n}{\Delta}$ and that $\ell > \frac{\Delta - 2}{\Delta} n$. The latter inequality implies that there is at least one leaf and, hence, at least one support vertex as well in $G$. In this case, the average number, $\frac{\ell}{s}$, of leaves  per support vertex satisfies $\frac{\ell}{s} > \Delta - 2$. Moreover, as $G$ is not isomorphic to a star, the maximum number of leaves adjacent to a support vertex is $\Delta-1$. Hence, there exists at least one support vertex which is adjacent to exactly $\Delta-1$ leaves. Therefore, we must have $G \cong G' \rhd_x S$, where $G'$ is also a tree, $x$ is a vertex of $G'$ and $S$ is a $\Delta$-star. Let $n' = |V(G')|$ and $\Delta'$ be the maximum degree of $G'$. We proceed further with the following claim.

\begin{unnumbered}{Claim~2.A}
If $n' \le 4$, then $\ID(G) \le ( \frac{\Delta - 1}{\Delta} ) n$.
\end{unnumbered}

\begin{proofofclaim}
Suppose that $n' \le 4$. If $n'=1$, then $G$ is a $\Delta$-star by itself and, hence, $G$ is isomorphic to a graph in $\mathcal{T}_\Delta$, a contradiction. If $n'=2$, then since $G'$ is connected, we have $G' \cong P_2$. Therefore, $G \cong P \rhd S$, where $P$ is a $2$-path and $S$ is a $\Delta$-star, and the desired upper bound follows by Lemma~\ref{lem_(m,Delta)-broom}(1). If $n'=3$, then since $G'$ is a tree, $G' \not \cong K_3$ and, thus, $G \in P \rhd S$, where $P$ is a $3$-path, and once again, by Lemma~\ref{lem_(m,Delta)-broom}(1), the desired result holds. Hence, we may assume that $n' = 4$.

Since $G'$ is a tree, $G'$ can be isomorphic to either a $4$-path or a $3$-star. Let us, therefore, assume that $G'$ is isomorphic to either a $4$-path $P$, say, or a $3$-star $S_1$, say. If $G' \cong P$, then by the fact that $G$ is not isomorphic to the graph $T_2$ in $\mathcal{T}_{3}$, we must have $G \cong P \rhd_x S$, where $x$ is a non-leaf vertex of $P$. In this case, by Lemma~\ref{lem_(m,Delta)-broom}(2), the result holds. If, however, $G' \cong S_1$, then $G \in S_1 \rhd_x S$, where $x$ is a vertex of $S_1$. If $x$ is a leaf of $S_1$ and $\Delta = 3$, then $G \cong T_1$ in the collection $\mathcal{T}_{3}$, a contradiction. Therefore, in the case that either $x$ is a non-leaf vertex of $S_1$ or that $\Delta \ge 4$, we are done by Lemma~\ref{lem_4-stars}.
\end{proofofclaim}

\medskip
By Claim~2.A, we may assume that $n' \ge 5$, for otherwise the desired result follows. Now, if $\Delta' = 2$, then $G'$ is a path of order at least $5$ and hence, we are done by Lemma~\ref{lem_(m,Delta)-broom}(1). Hence we may assume in what follows that $\Delta' \ge 3$.

\begin{unnumbered}{Claim~2.B}
If $G'$ is isomorphic to a graph in $\mathcal{T}_{\Delta'}$, then $\ID(G) \le ( \frac{\Delta - 1}{\Delta} ) n$.
\end{unnumbered}

\begin{proofofclaim}
Suppose that $G'$ is isomorphic to a graph in $\mathcal{T}_{\Delta'}$. If $G'$ is isomorphic to a star in $\mathcal{T}_{\Delta'}$, then $G' \cong S'$, where $S'$ is a $\Delta'$-star. Since $n' \ge 5$, we must have $\Delta' \ge 4$, implying that $G \cong S' \rhd S$, and the desired upper bound follows from Lemma~\ref{lem_4-stars}. If $G'$ is isomorphic to $T_2$, then $G = T_2 \rhd S$. However, as $G \not\cong T_3 \in \mathcal{T}_{3}$, the result follows by Lemma~\ref{lem_T2 rhd S}. If $G'$ is isomorphic to $T_3$, then in this case, $G = T_3 \rhd S$ and hence, the result holds by Lemma~\ref{lem_T3 rhd S}.

Hence, we may assume that $G'$ is isomorphic to an appended $3$-star in $\mathcal{T}_{3}$. Thus, $G$ is an appended star. If $\Delta \ge 4$, then we are done by Lemma~\ref{lem_4-stars}. Hence we may assume that $\Delta = 3$, implying that $G$ is an appended $3$-star itself. However, for $G$ not to be isomorphic to any graph in $\mathcal{T}_\Delta$, the tree $G$ must be an appended $3$-star of diameter at least~$8$ and in which case, by Lemma~\ref{lem_3-star diam 8}, we are done.
\end{proofofclaim}

\medskip
By Claim~2.B, we may assume that $G'$ is not isomorphic to any tree in $\mathcal{T}_{\Delta'}$, for otherwise the desired result follows. Since $n' \ge 5$ and $\Delta' \ge 3$, 
by the induction hypothesis, we have
\[
\ID(G') \le \left(  \frac{\Delta - 1}{\Delta} \right) n'.
\]
Thus, by Lemma~\ref{lem_ID(G) from ID(G')}, the result holds.
\end{proof}

\section{Tightness of the bound for trees}\label{sec:construct}

In this section, we consider some trees for which Conjecture~\ref{conj_G Delta_UB} is tight. Clearly, Conjecture~\ref{conj_G Delta_UB} is tight for every graph in $\mathcal{T}_\Delta$ (with $c=1/\Delta$).

Moreover, it is tight for some infinite families (with $c=0$), such as double stars with $2\Delta-2$ leaves, where by a \textit{double star} we mean a tree $G = S_1\rhd_u S_2$, where $S_1$ is a $(\Delta-1)$-star, $S_2$ is a $\Delta$-star and $u$ is universal vertex of $S_1$. The double star $G$ has $n=2\Delta$ vertices and the smallest identifying code $C$ of $G$ has $2\Delta-2$ code vertices. Thus, $|C|= ( \frac{\Delta-1}{\Delta} )n$.

Another class of trees for which the conjecture is tight when $\Delta=3$ is the $2$-corona of a path $P$ (see~\cite{FL22} or \cite[Section~1.3]{henning2013total}). We obtain the \emph{$2$-corona} of a graph $G$ by identifying each vertex $v$ of $G$ with a leaf of a $3$-vertex path $P_v$.

We also believe that more examples can be built using the appended stars and paths constructions from Section~\ref{sec:small_cases}.

Moreover, in the regime when both $n$ and $\Delta$ are large, we next show that there are trees which almost attain the conjectured bound. Notice that when $\Delta$ is large, the value $\frac{\Delta-1+\frac{1}{\Delta-2}}{\Delta+\frac{2}{\Delta-2}}$ is roughly $\frac{\Delta-1}{\Delta}$.

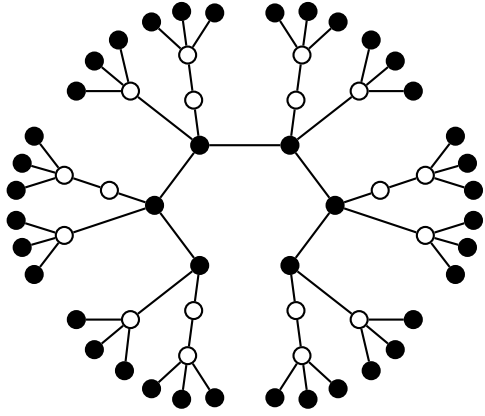
\begin{figure}[t!]
\centering
\begin{tikzpicture}[
blacknode/.style={circle, draw=black!, fill=black!, thick},
whitenode/.style={circle, draw=black!, fill=white!, thick},
scale=0.4]
\tiny
\node[blacknode] (2) at (-3,0) {}; 
\node[blacknode] (3) at (-1.5,2) {}; 
\node[blacknode] (4) at (-1.5,-2) {}; 
\node[blacknode] (5) at (1.5,2) {}; 
\node[blacknode] (6) at (1.5,-2) {}; 
\node[blacknode] (7) at (3,0) {}; 

\node[whitenode] (8) at (6,1) {}; 
\node[blacknode] (9) at (7,2.3) {};
\node[blacknode] (10) at (7.6,0.5) {};
\node[whitenode] (11) at (4.5,0.5) {};

\node[blacknode] (50) at (7.4,1.4) {};

\node[whitenode] (12) at (6,-1) {};
\node[blacknode] (13) at (7,-2.3) {};
\node[blacknode] (14) at (7.6,-0.5) {};

\node[blacknode] (51) at (7.4,-1.4) {};

\node[whitenode] (15) at (-6,1) {}; 
\node[blacknode] (16) at (-7,2.3) {};
\node[blacknode] (17) at (-7.6,0.5) {};
\node[whitenode] (18) at (-4.5,0.5) {};

\node[blacknode] (52) at (-7.4,1.4) {};

\node[whitenode] (19) at (-6,-1) {};
\node[blacknode] (20) at (-7,-2.3) {};
\node[blacknode] (21) at (-7.6,-0.5) {};

\node[blacknode] (53) at (-7.4,-1.4) {};

\node[whitenode] (22) at (-1.9,5) {}; 
\node[blacknode] (23) at (-1,6.4) {};
\node[blacknode] (24) at (-3.1,6.1) {};
\node[whitenode] (25) at (-1.7,3.5) {};

\node[blacknode] (54) at (-2.1,6.45) {};

\node[whitenode] (26) at (-3.8,3.8) {};
\node[blacknode] (27) at (-4.2,5.5) {};
\node[blacknode] (28) at (-5.6,3.8) {};

\node[blacknode] (55) at (-5,4.8) {};


\node[whitenode] (29) at (-1.9,-5) {}; 
\node[blacknode] (30) at (-1,-6.4) {};
\node[blacknode] (31) at (-3.1,-6.1) {};
\node[whitenode] (32) at (-1.7,-3.5) {};

\node[blacknode] (56) at (-2.1,-6.45) {};

\node[whitenode] (33) at (-3.8,-3.8) {};
\node[blacknode] (34) at (-4.,-5.5) {};
\node[blacknode] (35) at (-5.6,-3.8) {};

\node[blacknode] (57) at (-5,-4.8) {};


\node[whitenode] (36) at (1.9,5) {};
\node[blacknode] (37) at (1,6.4) {};
\node[blacknode] (38) at (3.1,6.1) {};
\node[whitenode] (39) at (1.7,3.5) {};

\node[blacknode] (58) at (2.1,6.45) {};

\node[whitenode] (40) at (3.8,3.8) {};
\node[blacknode] (41) at (4.2,5.5) {};
\node[blacknode] (42) at (5.6,3.8) {};

\node[blacknode] (59) at (5,4.8) {};


\node[whitenode] (43) at (1.9,-5) {}; 
\node[blacknode] (44) at (1,-6.4) {};
\node[blacknode] (45) at (3.1,-6.1) {};
\node[whitenode] (46) at (1.7,-3.5) {};

\node[blacknode] (60) at (2.1,-6.45) {};

\node[whitenode] (47) at (3.8,-3.8) {};
\node[blacknode] (48) at (4.2,-5.5) {};
\node[blacknode] (49) at (5.6,-3.8) {};

\node[blacknode] (61) at (5,-4.8) {};


\draw[-, thick, black!] (2) -- (3);
\draw[-, thick, black!] (3) -- (5);
\draw[-, thick, black!] (5) -- (7);
\draw[-, thick, black!] (7) -- (6);
\draw[-, thick, black!] (4) -- (2);


\draw[-, thick] (7) -- (11);
\draw[-, thick] (11) -- (8);
\draw[-, thick] (8) -- (9);
\draw[-, thick] (8) -- (10);
\draw[-, thick] (8) -- (50);

\draw[-, thick] (7) -- (12);
\draw[-, thick] (12) -- (13);
\draw[-, thick] (12) -- (14);
\draw[-, thick] (12) -- (51);


\draw[-, thick] (2) -- (18);
\draw[-, thick] (18) -- (15);
\draw[-, thick] (15) -- (16);
\draw[-, thick] (15) -- (17);
\draw[-, thick] (15) -- (52);

\draw[-, thick] (2) -- (19);
\draw[-, thick] (19) -- (20);
\draw[-, thick] (19) -- (21);
\draw[-, thick] (19) -- (53);


\draw[-, thick] (3) -- (25);
\draw[-, thick] (25) -- (22);
\draw[-, thick] (22) -- (23);
\draw[-, thick] (22) -- (24);
\draw[-, thick] (22) -- (54);

\draw[-, thick] (3) -- (26);
\draw[-, thick] (26) -- (27);
\draw[-, thick] (26) -- (28);
\draw[-, thick] (26) -- (55);


\draw[-, thick] (4) -- (32);
\draw[-, thick] (32) -- (29);
\draw[-, thick] (29) -- (30);
\draw[-, thick] (29) -- (31);
\draw[-, thick] (29) -- (56);

\draw[-, thick] (4) -- (33);
\draw[-, thick] (33) -- (34);
\draw[-, thick] (33) -- (35);
\draw[-, thick] (33) -- (57);


\draw[-, thick] (5) -- (39);
\draw[-, thick] (39) -- (36);
\draw[-, thick] (36) -- (37);
\draw[-, thick] (36) -- (38);
\draw[-, thick] (36) -- (58);

\draw[-, thick] (5) -- (40);
\draw[-, thick] (40) -- (41);
\draw[-, thick] (40) -- (42);
\draw[-, thick] (40) -- (59);


\draw[-, thick] (6) -- (46);
\draw[-, thick] (46) -- (43);
\draw[-, thick] (43) -- (44);
\draw[-, thick] (43) -- (45);
\draw[-, thick] (43) -- (60);

\draw[-, thick] (6) -- (47);
\draw[-, thick] (47) -- (48);
\draw[-, thick] (47) -- (49);
\draw[-, thick] (47) -- (61);

\end{tikzpicture}
\caption{\small Tree $T_{t,\Delta}$ as in Proposition~\ref{prop:BigConstruction} with $t=3$ and $\Delta = 4$. The set of black vertices constitutes an identifying code of $T_{t,\Delta}$.}
\label{fig_tightness}
\end{figure}

\begin{proposition}\label{prop:BigConstruction}
Let $P_{t}$ be a path of order $t\ge 3$ and let $\Delta\ge 4$ be an integer. Let an intermediate tree of order~$n$ be formed by appending onto every vertex of the path $\Delta-2\ge 2$ copies of $\Delta$-stars. Thereafter, for each vertex $p_i$ of the path, subdivide a single edge between $p_i$ and an adjacent support vertex of the intermediate tree. If $T_{t,\Delta}$ denotes the resulting tree (see Figure~\ref{fig_tightness} for an example with $t=6$ and $\Delta = 4$), then
\[
\ID(T_{t,\Delta}) =
\left( \frac{\Delta-1+\frac{1}{\Delta-2}}{\Delta+\frac{2}{\Delta-2}} \right) n > \left( \frac{\Delta-1}{\Delta} \right) n  - \frac{n}{\Delta^2} .
\]
\end{proposition}
\begin{proof}
Let $T=T_{t,\Delta}$. For each $j\in[t]$, let $p_j$ be a vertex in the path $P_{t}$ and let $S^j_i$ for $i \in [\Delta-2]$ be $\Delta$-stars appended onto $p_j$. Moreover, let $u^j_1$ be the universal vertex of $S^j_1$ and let the edge $u_1^jp_j$ of the intermediate tree be subdivided by the vertex $v_1^j$ of $T$. We then have $n=t(\Delta-2)\Delta+2\cdot t$ and the maximum degree of $T$ is $\Delta$. Observe that any identifying code of $T$ requires at least $\Delta-1$ codewords from each star $S^j_i$. Furthermore, we need an additional codeword to dominate $v^j_1$ (if the center of $S^j_1$ is not chosen yet) or to separate $v_1^j$ from a leaf of $S_1^j$ (if the center of $S_1^j$ is chosen but some leaf of $S_1^j$ is not). Thus, there are at least $(\Delta-2)(\Delta-1)+1$ codewords among the vertices in $V_j(P_{t}) = \{p_1,v_1^j\} \cup V(S_1^j) \cup V(S_2^j) \cup \cdots \cup V(S_{\Delta-2}^j)$. Moreover, we have $|V_j(P_{t})| = (\Delta-2)\Delta+2$. Thus,
\[
\frac{\ID(T)}{n} \ge \frac{(\Delta-2)(\Delta-1)+1}{(\Delta-2)\Delta+2}=\frac{\Delta-1+\frac{1}{\Delta-2}}{\Delta+\frac{2}{\Delta-2}}.
\]

Moreover, we notice that we may choose as an identifying code of this size, all the leaves and all the path vertices (as illustrated in  Figure~\ref{fig_tightness} by the set of black vertices). Therefore,
\[
\frac{\ID(T)}{n} = \frac{\Delta-1+\frac{1}{\Delta-2}}{\Delta+\frac{2}{\Delta-2}}.
\]

We note that
\begin{align*}
\frac{\Delta-1+\frac{1}{\Delta-2}}{\Delta+\frac{2}{\Delta-2}}
     =& \frac{1}{\Delta} \left( \frac{\Delta^2-\Delta+\frac{\Delta}{\Delta-2}}{\Delta+\frac{2}{\Delta-2}} \right) \2 \\
     =& \frac{1}{\Delta} \left(   \frac{\Delta^3-3\Delta^2+3\Delta}{\Delta^2-2\Delta+2} \right) \2 \\
     =&\frac{\Delta-1}{\Delta}-\frac{\Delta-2}{\Delta^3-2\Delta^2+2\Delta} \2\\
     >&\frac{\Delta-1}{\Delta}-\frac{\Delta-2}{\Delta^2(\Delta-2)} \2\\
     =&\frac{\Delta-1}{\Delta}-\frac{1}{\Delta^2}.
\end{align*}
Hence, the claimed inequality holds.
\end{proof}

In the previous proposition, we could obtain a slightly better result by adding an additional star for each end of the path. However, that improvement would be of local nature and hence, would improve the result only by a constant while complicating the proof. Note that we may modify the tree $T_{t,\Delta}$ by adding an edge between the two endpoints of the underlying path $P_t$. The resulting graph has the same identification number as the tree $T_{t,\Delta}$.

\medskip
Previously, in \cite{BCHL2004}, an exact value for the identification number of complete $q$-ary trees was given. Later in \cite{foucaud2012size}, the authors showed that in terms of maximum degree, the complete $(\Delta-1)$-ary tree $T_\Delta$ has identification number
\[
\ID(T_\Delta)=\left\lceil \left( \frac{\Delta-2+\frac{1}{\Delta}}{\Delta-1+\frac{1}{\Delta}} \right) n\right\rceil.
\]
When we compare this value to the one in Proposition~\ref{prop:BigConstruction}, we observe that it has slightly smaller multiplier in front of $n$.

\medskip

An observant reader may notice that the root of the tree $T_\Delta$ has degree of only $\Delta-1$. Hence, we might gain slight improvements by adding some new subtree to it. Let $r$ be the root of tree $T_\Delta$ and let $T'$ be another tree which we join with an edge to $r$. Let $C_\Delta$ (resp., $C'$) be an identifying code in $T_\Delta$ (resp. $T'$). Then, $C=C_\Delta\cup C'\cup\{r\}$ is an identifying code in the combined tree. Furthermore, we may interpret the size of the identifying code also as the \textit{density} of the identifying code in the graph. In particular, we can observe that the density of the code in the combined tree is increased by a meaningful amount if and only if tree $T'$ is large and a minimum identifying code of $T'$ has a density larger than that of a minimum identifying code in $T_\Delta$. In other words, the tree $T'$ alone is a better example for a tree with large identifying code than the tree $T_\Delta$. Hence, tree $T_\Delta$ does not help in finding a larger example than we have found in Proposition~\ref{prop:BigConstruction}.

\medskip
Another observation we make is that when $n$ is large and $\Delta=4$, Proposition~\ref{prop:BigConstruction} gives $\ID(T_{t,4})=\frac{7}{10}n$. We do not know any large triangle-free constructions with significant improvements for this case.

\section{Conclusion}\label{sec:conclu}

We have made significant progress towards Conjecture~\ref{conj_G Delta_UB} by proving it for all trees; moreover we have precisely characterized the trees that need $c>0$ in the conjectured bound. It is interesting that for every fixed $\Delta \ge 3$, there is only a finite list of such trees. We close with the following list of open questions and problems.

\begin{question}
\label{question1}
For every fixed $\Delta \ge 3$, if $G$ is a connected identifiable graph of order $n$ and of maximum degree $\Delta$, then is it true that
\[
\ID(G) \le \left( \frac{\Delta - 1}{\Delta} \right) n ,
\]
except for a finite family of graphs?
\end{question}

\begin{problem}
\label{problem1}
Characterize the trees $T$ of order~$n$ with maximum degree~$\Delta$ satisfying
\[
\ID(T) = \left( \frac{\Delta - 1}{\Delta} \right) n.
\]
\end{problem}

\begin{problem}
\label{problem2}
Characterize the trees $T$ of order~$n \ge 3$ satisfying $\ID(T) = n-\gamma(T)$, that is, characterize the trees $T$ that achieve equality on the upper bound in Theorem~\ref{the_domBound}.
\end{problem}

\begin{problem}
\label{problem3}
Determine classes of graphs $G$ of order~$n \ge 3$ satisfying $\ID(T) \le n-\gamma(T)$, that is, does Theorem~\ref{the_domBound} hold for a larger class of graphs (for example some subclass of bipartite graphs containing all trees)?
\end{problem}

As mentioned before, in the companion paper~\cite{PaperPart2}, we use the findings from the current paper, to prove Conjecture~\ref{conj_G Delta_UB} for all triangle-free graphs (with the same list of graphs requiring $c>0$ when $\Delta\ge 3$). To do so, we will use Theorem~\ref{the_main} as the starting point of the proof. For $\Delta\ge 3$, the extremal triangle-free graphs (requiring $c>0$) are actually the same extremal trees characterized in the current paper.

\bibliographystyle{abbrv}
\bibliography{main_trees_copy-edited}

\end{document}